\documentclass{article}

\usepackage{amssymb}
\usepackage{natbib}
\usepackage{bm}
\usepackage{subfigure}
\usepackage{graphicx}
\usepackage[hmargin=3cm,vmargin=3.5cm]{geometry}
\usepackage{color}
\usepackage{amsthm}
\usepackage{amsmath}
\usepackage{booktabs}
\usepackage{multirow}

\newcommand{\vbar}{\,|\, }

\newtheorem{theorem}{Theorem}

\newtheorem{lemma}{Lemma}
\newtheorem{proposition}{Proposition}

\newcommand{\Prb}{\mathbb{P}}
\newcommand{\E}[2]{ {\mathbb E}_{#1} \left[ #2 \right] }
\begin{document}

\title{Log-Scale Shrinkage Priors and Adaptive Bayesian Global-Local Shrinkage Estimation}
\author{Daniel F. Schmidt\footnote{Email: {\tt daniel.schmidt@monash.edu}} \\ Faculty of Information Technology \\ Monash University \\ Australia \and Enes Makalic\footnote{Email: {\tt emakalic@unimelb.edu.au}} \\ School of Population and Global Health \\ The University of Melbourne \\ Australia} 
\maketitle 

\begin{abstract}

Global-local shrinkage hierarchies are an important innovation in Bayesian estimation. We propose the use of log-scale distributions as a novel basis for generating familes of prior distributions for local shrinkage hyperparameters. By varying the scale parameter one may vary the degree to which the prior distribution promotes sparsity in the coefficient estimates. By examining the class of distributions over the logarithm of the local shrinkage parameter that have log-linear, or sub-log-linear tails, we show that many standard prior distributions for local shrinkage parameters can be unified in terms of the tail behaviour and concentration properties of their corresponding marginal distributions over the coefficients $\beta_j$. We derive upper bounds on the rate of concentration around $|\beta_j|=0$, and the tail decay as $|\beta_j| \to \infty$, achievable by this wide class of prior distributions. 

We then propose a new type of ultra-heavy tailed prior, called the log-$t$ prior with the property that, irrespective of the choice of associated scale parameter, the marginal distribution always diverges at $\beta_j = 0$, and always possesses super-Cauchy tails. We develop results demonstrating when prior distributions with (sub)-log-linear tails attain Kullback--Leibler super-efficiency and prove that the log-$t$ prior distribution is always super-efficient. We show that the log-$t$ prior is less sensitive to misspecification of the global shrinkage parameter than the horseshoe or lasso priors. By incorporating the scale parameter of the log-scale prior distributions into the Bayesian hierarchy we derive novel adaptive shrinkage procedures. Simulations show that the adaptive log-$t$ procedure appears to always perform well, irrespective of the level of sparsity or signal-to-noise ratio of the underlying model.\\

\end{abstract}

\newpage

\def\T{{ \mathrm{\scriptscriptstyle T} }}

\section{Introduction}
\label{sec:intro:properties}
The multiple means problem has been studied extensively since the introduction of the first shrinkage estimator by \cite{JamesStein61}. In the multiple means problem we observe a vector of $n$ samples ${\bf y} = (y_1,\ldots,y_n)$ from the model
\begin{equation}
	\label{eq:multiple:means}
	y_i \vbar \beta_i, \sigma \sim N(\beta_i, \sigma^2)
\end{equation}
and are required to estimate the unknown coefficient (mean) vector $\bm{\beta} = (\beta_1,\ldots,\beta_n)$. Part of the appeal of the multiple means problem is that it serves as an analogue for the more complex general linear model, while being substantially more amenable to analysis. The original work of James and Stein~\citep{JamesStein61} showed that proportional shrinkage can be used to construct estimators that dominate least-squares in terms of squared-error risk for $n \geq 3$, and these results were later extended to Bayesian proportional shrinkage estimators~\citep{Strawderman71,Zellner86a}, and more generally to Bayesian ridge regression. More recently there has been a focus on the sparse setting in which the majority of the entries of $\bm{\beta}$ are exactly zero. In this setting the use of methods that promote sparsity, such as the lasso~\citep{Tibshirani96} can lead to substantial improvements in estimation risk over conventional shrinkage estimators. An important contribution to the Bayesian regression literature was the introduction of the general global-local shrinkage hierarchy proposed in \cite{PolsonScott10}, which models the coefficients by
\begin{eqnarray}
	\label{eq:Global}
	\beta_j \vbar \lambda_j, \tau &\sim& N(0, \lambda_j^2 \tau^2 \sigma^2), \\
	\lambda_j &\sim& p(\lambda_j) d\lambda_j, \\
	     \tau &\sim& p(\tau) d\tau	.
	\label{eq:Local}
\end{eqnarray}
In the global-local shrinkage hierarchy the local shrinkage hyperparameters $\lambda_j$ control the degree of shrinkage applied to individual coefficients, while the global shrinkage hyperparameter $\tau$ controls the overall degree of shrinkage applied globally to all coefficients. The choice of the prior distributions $p(\lambda_j)$ controls the behaviour of the resulting shrinkage estimators; for example, if $p(\lambda_j)$ is a Dirichlet point-mass the hierarchy reduces to the standard ridge regression prior. Well known techniques that fall under the umbrella of global-local shrinkage priors include the Bayesian lasso~\citep{ParkCasella08}, the normal-gamma prior~\citep{GriffinBrown10}, the horseshoe~\citep{CarvalhoPolson10}, the horseshoe+~\citep{BhadraEtAl15}, the beta mixture of Gaussians~\citep{ArmaganEtAl11b} and the R2-D2 prior~\citep{ZhangEtAl16}. As there exists such a wide array of potential global-local shrinkage priors, it is valuable to determine general properties that may or may not be beneficial for estimating $\bm{\beta}$ in the multiple means problem. \cite{CarvalhoPolson10} proposed two such properties that a good sparsity inducing prior should possess:
\begin{itemize}
	\item {\em Property I:} The prior $p(\lambda_j)$ should concentrate sufficient probability mass near $\lambda_j=0$ to ensure that the marginal distribution $p(\beta_j)$ has an infinitely tall spike at the origin; that is, $p(\beta_j) \to \infty$ as $|\beta_j| \to 0$;
	
	\item {\em Property II:} The prior $p(\lambda_j)$ should decay sufficiently slowly as $\lambda_j \to \infty$ so that the marginal distribution $p(\beta_j)$ has Cauchy, or super-Cauchy tails, to ensure that
	\[
		\E{}{\beta_j \vbar {\bf y}} = y_j + o(1)
	\]
	as $|y_j| \to \infty$.
\end{itemize}
The first property ensures that the Bayes estimation risk of a procedure using such a prior will be low when the underlying coefficient vector is sparse. The second property ensures that very large effects are not over-shrunk by the resulting procedure. A number of shrinkage priors, including the horseshoe and horseshoe+, satisfy both of these properties, and also possess a number of other favourable theoretical properties~(see for example, \cite{VanDerPasEtAl17a}). More classical shrinkage priors, such as the Bayesian lasso, and Bayesian ridge do not satisfy either Property I or II; despite this, it is not difficult to propose configurations of the underlying coefficients $\bm{\beta}$ for which sparse estimation methods possess much greater squared-error risk than the humble Bayesian ridge. Generally the coefficient configurations that are problematic for sparse estimators are dense in the sense that a large number of the coefficients are non-zero. While it can be argued that sparsity inducing priors are inappropriate for such settings it is difficult to know with certainty whether a problem is best modelled {\em a priori} as sparse or dense, particularly when dealing with complex natural phenomena.

Statistical genetics problems have been an important application, and source of inspiration, for much of the recent work on sparsity inducing priors for high dimensional regression models. A key assumption driving much of this original work was that only a small number of the many genomic variants that exist are associated with any given disease. However, there is substantial evidence that there can be large numbers of variants associated with diseases~\citep{BoyleEtAl17a}, though the levels of association, and therefore the signal-to-noise ratio, will be low. Thus, there is a potential need for shrinkage priors that can adapt to the sparsity of the underlying coefficient vector, and bridge the gap between extreme sparsity inducing behaviour, and proportional shrinkage type priors, as the problem demands. Many of the existing sparsity inducing shrinkage priors have shape hyperparameters that can be tuned to adjust the degree to which they expect the underlying coefficient vector to be sparse, and some work has examined the tuning, and estimation, of these type of shape parameters~(see for example, \cite{BhattacharyaEtAl15}, \cite{GriffinBrown16} and \cite{HuberFeldkircher17a}). However, in general, the shape parameters are complex to sample, their interpretation is potentially difficult, and the tuning has often been focused on increasing the degree of {\em a priori} expected sparsity above and beyond that of a method such as the horseshoe. 

In this paper we propose a framework for specifying prior distribution for local shrinkage hyperparameters based on log-scale distributions. These log-scale distributions have a single scale hyperparameter that can be used to vary the log-scale shrinkage priors from highly sparsity promoting through to almost ridge regression-like in behaviour, and are straightforward to integrate into an MCMC procedure. Apart from providing a degree of robustness to our prior assumptions regarding sparsity of the underlying coefficients, we believe that these types of adaptive shrinkage priors will have particular application in regression problems in which variables can be formed into logical groupings. In such situations, it is highly concievable that the coefficient vectors associated with some groups will be dense while the coefficient vectors associated with other groups could be very sparse. Some obvious examples of this include additive models based around polynomial expansions in which one of the input variables may be related to the target through a discontinuous non-linearity, or statistical genomics, in which some genes or pathways may have large numbers of associations with disease, such as the HLA region~\citep{KennedyEtAl17a} while other genes may have only one or two strongly associated variants.





\subsection{Our Contribution}

In this paper we show that viewing prior distributions in terms of the logarithm of the local shrinkage parameters, $\xi_j = \log \lambda_j$, has several distinct advantages. Our work was motivated by the observation that those shrinkage priors which strongly promoted sparsity spread their probability mass more thinly across the $\xi_j$ space. By viewing the standard prior distributions in $\xi_j$ space and introducing a scale parameter $\psi$ it becomes possible to vary the degree to which a prior promotes sparsity in the coefficient estimates, all the way from the simple proportional shrinkage ridge regression model up to extremely heavy tailed distributions. We call these {\em log-scale} prior distributions.


Using this approach we show that many of standard local shrinkage parameter prior distributions can be unified in terms of tail behaviour and concentration properties of the resulting marginal distribution over $\beta_j$. In particular, we consider the class of distributions over $\xi_j$ that have log-linear, or sub-log-linear tails. We derive upper bounds on the rate of concentration around $|\beta_j|=0$, and tail decay as $|\beta_j| \to \infty$, achievable by this class of prior distributions. Further, we show that by introduction of a scale parameter, all of the common prior distributions can be made to behave equivalently to each other, irrespective of the specific shape parameters they may possess, in the sense that for a sufficiently small choice of scale $\psi$ the induced marginal distribution can be made to lose Properties I and II. We then propose a new class of ultra-heavy tailed priors, called the log-$t$ priors, which exhibit the property that irrespective of the choice of $\psi$, the induced marginal distribution over $\beta_j$ never loses Properties I and II. We further develop results demonstrating when prior distributions with (sub)-log-linear tails will attain Kullback--Leibler super-efficiency, as well as showing that the log-$t$ prior distribution is always super-efficient in this sense. Additionally, we show that the log-$t$ prior is less sensitive to misspecification of $\tau$ than the horseshoe (and related distributions with (sub)-log-linear tails).

Finally, we utilise the simple interpretation of $\psi$ as the scale for the $\xi_j$ hyper-parameters to derive an adaptive shrinkage procedure. We incorporate $\psi$ into the full Bayesian hierarchy and use this to estimate the degree of sparsity in the data generating model. This yields prior distributions that are able to vary from highly sparsity promoting through to ridge-like, depending on the configuration of the true regression coefficients. However, by using the log-$t$ prior distribution the resulting prior distribution over the coefficients $\beta_j$ never loses Properties I and II no matter how much mass is concentrated near $\xi_j = 0$.

\section{Log-scale Hyperprior Distributions}
\label{sec:LogScale}

Let $f(\xi_j)$ be a unimodal distribution over $\mathbb{R}$. If $\xi_j \sim f(\xi_j) d \xi_j$, then the translated and scaled random variate $\xi_j^\prime = \sigma \xi_j + \eta$ follows the probability distribution:
%
\begin{equation}
	\label{eq:location:scale}
	p(\xi_j^\prime \vbar \eta, \psi) = \left( \frac{1}{\psi} \right) f \left( \frac{\xi_j^\prime - \eta}{\psi} \right).
\end{equation}
Distributions of the form (\ref{eq:location:scale}) are known as location-scale distributions, in which $\eta \in \mathbb{R}$ is the location parameter and $\psi \in \mathbb{R}_+$ is the scale parameter. Let $\xi_j = \log \lambda_j$ be the natural logarithm of local shrinkage parameter $\lambda_j$ in the global-local shrinkage hierarchy (\ref{eq:Global})--(\ref{eq:Local}). The primary motivation for studying distributions over $\xi_j$ is the fact that, if $\xi_j$ follows a location-scale distribution of the form (\ref{eq:location:scale}), then:
\begin{enumerate}
	\item location transformations of $\xi_j$ induce scale transformations on $\lambda_j$;
	\item scale transformations of $\xi_j$ induce power-transformations on $\lambda_j$. 
\end{enumerate}
The first fact is of less interest, as scale transformations of $\lambda_j$ are generally taken care of by the presence of the global shrinkage hyperparameter $\tau^2$ in the standard global-local shrinkage prior hierarchy. The second fact is far more interesting, as it reveals a simple way in which we can control the behaviour of the prior distribution $p(\lambda_j)$ when $\lambda_j \to 0$ and $\lambda_j \to \infty$. If we further restrict attention to the class of log-location-scale priors in which $f(\xi_j)$ is symmetric around $\xi_j=\eta$ the prior distribution $p(\lambda_j)$ is symmetric around $\lambda_j=e^\eta$ in the following sense:
	\begin{equation}
		\label{eq:Xi:symmetry}
		\Prb \left( \lambda_j \in (e^\eta/k, e^\eta) \right) = \Prb \left( \lambda_j \in (e^\eta, k e^\eta) \right), \; \; k>0.
	\end{equation}
A property of the global-local shrinkage hierarchy for the multiple means problem is that the posterior mean of the coefficients $\beta_j$ can be written as
\[
	\E{}{\beta_j \vbar {\bf y}} = \left( 1-\E{}{\kappa_j \vbar {\bf y}} \right) \, y_j
\]
where $\kappa_j$ is the degree of shrinkage towards zero being applied to coefficient $\beta_j$. Given $\lambda_j$, the corresponding degree of shrinkage is $\kappa_j = 1/(1+\lambda_j^2)$; when $\kappa_j$ is close to zero very little shrinkage is performed, and when $\kappa_j$ is close to one, the corresponding coefficient is almost entirely shrunk to zero. This interpretation motivated the original horseshoe prior distribution, which placed a horseshoe-shape prior over $\kappa_j$ to promote either aggressive shrinkage or little shrinkage of the coefficients. The quantity $\E{}{\kappa_j \vbar {\bf y}}$ can be interpreted as the degree of evidence in the data to support $\beta_j \neq 0$ against $\beta_j = 0$, and the thresholding rule $\E{}{\kappa_j \vbar {\bf y}} < 1/2$ is frequently used as a variable selection criterion~\citep{CarvalhoPolson10,TangEtAl16}.

Placing a log-scale distribution over $\xi_j$ that is symmetric around $\eta$ implies a distribution $p(\kappa_j)$ over $\kappa_j$ with a median at $\kappa_j=1/(1+e^{2\eta})$. In the particular case that $\eta=0$, property (\ref{eq:Xi:symmetry}) implies that the resulting distribution models the prior belief that coefficients are just as likely to be shrunken towards zero as they are to be left untouched, and that {\em a priori}, a variable has a marginal prior probability of being selected of $1/2$. These properties, coupled with the fact that scale transformations of $\xi_j$ result in power-transformations of $\lambda_j$ suggest that specification of symmetric priors over $\xi_j$ may offer a fruitful approach to generate novel, adjustable priors for local shrinkage parameters that imply reasonable prior beliefs about the model coefficients.
	
%

\subsection{Behaviour of Standard Shrinkage Priors in $\xi_j$ Space}
\label{sec:StandardPriors}

It is of interest to examine the prior distributions $p(\xi_j)$ implied by a number of standard shrinkage prior distributions for $\lambda_j$. The Bayesian lasso (double exponential) \cite{ParkCasella08} prior distribution over $\beta_j$ induces an exponential distribution over $\lambda_j^2$, and a distribution of the form
\begin{equation}
	\label{eq:BayesLassoXi}
	p_{\rm DE}(\xi_j) = 2 e^{2 \xi_j - e^{2 \xi_j}}
\end{equation}
over $\xi_j$. This is an asymmetric distribution over $\xi_j$, with the left-hand tail (that controls shrinkage less than $1/2$) being much heavier than the right-hand tail (which controls shrinkage greater than $1/2$). This interquartile interval (first and third quartiles) for this distribution is approximately $(-0.623, \; 0.163)$. Positive values of $\xi_j$ induce little shrinkage on coefficients, which is desirable for modeling very large effects. The skew towards negative values of $\xi_j$ exhibited by the Bayesian lasso demonstrates why it introduces bias in estimation when the underlying model coefficients are large. In terms of the coefficient of shrinkage, $1-\kappa_j$, where $\kappa_j = 1/(1+e^{2\xi_j})$, the interquartile interval of the prior induced by (\ref{eq:BayesLassoXi}) is approximately $(0.223, \; 0.581)$, respectively, demonstrating a clear preference towards shrinkage below $1/2$. 

The horseshoe prior~\citep{CarvalhoPolson10} is often considered a default choice for sparse regression problems. The horseshoe prior places a standard half-Cauchy distribution over $\lambda_j$, which is known to induce an unstandardised unit hyperbolic secant distribution over $\xi_j$, with probability distribution given by
\begin{equation}
	\label{eq:HS:xi}
	p_{\rm HS}(\xi_j) = \left( \frac{1}{\pi} \right) {\rm sech}(\xi_j),
\end{equation}
where ${\rm sech}(\cdot)$ denotes the hyperbolic secant function. This distribution is symmetric around $\xi_j = 0$, and has an interquartile range of approximately $(-0.881, \; 0.881)$ for $\xi_j$, and $(0.15, \; 0.85)$ for $1-\kappa_j$. In contrast to the Bayesian lasso, the horseshoe prior clearly spreads its probability mass more thinly across the $\xi_j$ space, encodes a much wider range of a prior plausible shrinkage for the coefficients and is symmetric around $\xi_j=0$. The horseshoe+ (HS+) prior distribution~\cite{BhadraEtAl15} models $\lambda_j$ as the product of two half-Cauchy distributions, which leads to a prior distribution on $\xi_j$ of the form
\begin{equation}
	\label{eq:HSplus:xi}
	p_{\rm HS+}(\xi_j) = \left( \frac{2 \xi_j}{\pi^2} \right) {\rm csch}(\xi_j)
\end{equation}
where ${\rm csch}(\cdot)$ is the hyperbolic cosecant function. This distribution is also symmetric around $\xi_j = 0$, which is straightforward to verify from the fact that $\xi_j$ is modelled as the sum of two hyperbolic secant random variables, which are themselves symmetric. The HS+ prior has an interquartile range of approximately $(-1.33, \;  1.33)$ for $\xi_j$, which translates to an interquartile interval of $(0.062, \; 0.938)$ on $1-\kappa_j$. The horseshoe+ prior more strongly promotes sparsity in the estimates of $\bm{\beta}$, and this is evident by the fact that more probability mass is concentrated in a larger region of $1-\kappa_j$, than for either the Bayesian lasso or the horseshoe.
%
%
More generally, the beta prime class of hyperpriors, which begin by modelling the shrinkage factor as $\kappa_j \sim {\rm Be}(a,b)$~\citep{ArmaganEtAl11b}, imply a distribution of the form
\begin{equation}
	\label{eq:BP:tj}
	p_{\rm BP}(\xi_j \vbar a,b) = \frac{2 \, \Gamma(a+b)}{\Gamma(a) \Gamma(b)} \frac{e^{2 a \xi_j}}{\left( e^{2 \xi_j}+1 \right)^{a+b}}
\end{equation}
over $\xi_j$. The density (\ref{eq:BP:tj}) be identified as a $z$-distribution with zero mean, a scale of $1/2$ and shape parameters $a$ and $b$. This class of distributions is symmetric if and only if $a=b$, and generalizes a number of standard shrinkage hyperpriors. For example, the standard horseshoe is recovered by taking $a=b=1/2$, while the Strawderman-Berger prior is recovered by taking $a=1/2$ and $b=1$, which is asymmetric. The negative-exponential-gamma prior is found by taking $a=1$ and $b=c-2$, with $c>0$, which is asymmetric for all $c \neq 3$. In all these cases the hyperparameter $a$ controls the behaviour of the prior on $\lambda_j$ as $\lambda_j \to 0$, and the $b$ hyperparameter controls the behaviour of the tail of the prior on $\lambda_j$ as $\lambda_j \to \infty$. By adjusting the $a$ and $b$ hyperparameters, the prior mass can be controlled to be spread more or less densely across $\xi_j$ space, controlling the degree to which the prior induces sparsity on the estimated coefficients. 

For example, when $a=b=1/2$ the interquartile interval for $\xi_j$ is identical to the interval obtained for the horseshoe, while for $a = 1/4$ and $b = 1/4$ the interquartile range on $\xi_j$ expands to $(-1.53, \; 1.53)$ which is wider than interquartile range for the horseshoe+ prior. Taking $a>1$, $b>1$ leads to a concentration of prior mass near $\xi_j=0$, which can be used to approximate ridge regression. However, a potentially unwanted side effect of this is that when $a>1, b>1$, the marginal prior distribution for $\beta_j$ loses both Properties I and II. The effect of $a$ and $b$ on the prior distribution has been used to attempt to adaptively estimate the degree of sparsity required from the data (for example, \cite{GriffinBrown16} and \cite{HuberFeldkircher17a}). However, the functional form of the prior distribution, and the way in which the hyperparameters $a$ and $b$ control the prior, has the consequence that both the interpretation of the the hyperparameters, and the practical implementation of efficient sampling algorithms for them, is difficult.


\subsection{Log-Scale Priors as Shrinkage Priors on $\xi_j$}
\label{sec:log:scale:shrinkage:xi}

The log-scale interpretation of standard shrinkage priors offers an alternative way of understanding the way in which both the tails, and behaviour near the origin, of a prior distribution for $\lambda_j$ models prior beliefs regarding sparsity of $\bm{\beta}$, and the type of shrinkage behaviour introduced by the prior distribution. The standard prior distributions discussed in Section \ref{sec:StandardPriors} all induce unimodal distributions on $\xi_j$ which tail off to zero as $|\xi_j| \to \infty$. They differ in how thinly they spread their prior probability across the $\xi_j$ space. In the standard global-local shrinkage hierarchy (\ref{eq:Global})--(\ref{eq:Local}), the prior distribution for $\beta_j$ is
\[
	\beta_j \sim N(0, \lambda_j^2 \tau^2 \sigma^2).
\]
If we assume that $\sigma=1$, we see that conditional on $\tau$, the local shrinkage parameter $\lambda_j$ is modelled as a random variable scaled by $\tau$. This scale transformation of $\lambda_j$ induces a location transformation on $\xi_j = \log \lambda_j$; i.e., if $f(\xi_j)$ is the density implied by the prior for $\lambda_j$ over $\xi_j$, and $\xi_j^\prime = \log \lambda_j \tau$, then
\begin{equation}
	\label{eq:location:only}
	p(\xi_j^\prime \vbar \tau) = f(\xi_j^\prime - \log \tau).
\end{equation}
The global scale parameter determines the location of the prior over $\xi_j$. The standard shrinkage priors on $\lambda_j$ can therefore be viewed as shrinking the $\xi_j$ hyperparameters towards $\log \tau$, with the more sparsity promoting prior distributions resulting in less shrinkage of the $\xi_j$ hyperparameters. Clearly (\ref{eq:location:only}) is of the form (\ref{eq:location:scale}) with $\eta=\log \tau$ and $\psi=1$. A natural generalization is then to allow the scale $\psi$ to take on an arbitrary value. The scale $\psi$ of the log-scale prior (\ref{eq:location:scale}) can be interpreted as modelling the {\em a priori} plausible range of $\xi_j$ values around the location $\eta = \log \tau$. The smaller the scale parameter, the more prior probability is concentrated around $\xi_j = \eta$, and the less variability is implied in the values of $\xi_j$, with the result that most of the shrinkage coefficients will be concentrated around $\kappa = 1/(1+\tau^2)$. In the limiting case that $\psi \to 0$, the prior (\ref{eq:location:scale}) concentrates all of its mass at $\xi_j = \eta$, allowing for no variation in shrinkage between coefficients, and the prior hierarchy reduces to the Bayesian ridge. In contrast, the larger the scale parameter $\psi$ becomes, the more a prior variability in the $\xi_j$ hyperparameters is implied, with the caveat that less prior mass is placed around the neighbourhood of any particular $\xi_j$. In the limiting case that $\psi \to \infty$, we recover the (improper) normal-Jeffreys prior $1/\lambda_j$, which is a uniform distribution over $\xi_j$. 

This interpretation motivates us to propose the introduction of a scale parameter to the prior distributions over $\xi_j$ space as a method to provide a hyperparameter that can be used to control the amount by which a prior promotes sparsity. Practically, this type of scale hyperparameter is easier to deal with than shape hyperparameters that control the tail behaviour of priors such as the beta prime prior (\ref{eq:BP:tj}), both in terms of interpretation as well as implementation with a sampling hierarchy. They also provide a unified form of hyperparameter that controls the behaviour of a shrinkage prior in the same, standard way, irrespective of the initial prior distribution $f(\xi_j$) that we start with.

\section{Three Log-Scale Prior Distributions}

In this section we examine three potential choices of log-scale prior distribution for global-local shrinkage hierarchies of the form (\ref{eq:Global})--(\ref{eq:Local}). The first is the log-hyperbolic secant prior, which is itself a generalization of the regular horseshoe prior distribution. 
The second distribution we consider is the asymmetric log-Laplace prior distribution, which has the advantage of being amenable to analysis, while also exhibiting the same tail properties as the log-hyperbolic secant prior. Furthermore, the log-Laplace prior distribution can be used to derive upper-bounds on the concentration and tail behaviours of a large class of prior distributions, which includes most of the common shrinkage priors. The final distribution we consider is the log-$t$, which is formed by modelling the $\xi_j$ hyperparameters using the Student-$t$ distribution. The resulting density appears to be part of an entirely new class of prior distributions, and exhibits a special property that is, to the authors knowledge, not shared by any other known shrinkage prior.

\subsection{Log-Hyperbolic Secant Prior}
\label{sec:LogHyperbolicSecant}

The horseshoe prior is generally considered a default choice of prior distribution for the local shrinkage parameters $\lambda_j$, and therefore forms a suitable starting point for generalisation through the introduction of scale parameter on the $\xi_j$ space. Our starting point is the density (\ref{eq:HS:xi}), which after introduction of a location parameter $\eta$ and scale parameter $\psi$ becomes
\begin{equation}
	\label{eq:hyp:sech}
	p_{\rm SECH}(\xi_j \vbar \psi) = \left( \frac{1}{\pi \psi} \right) {\rm sech}\left( \frac{\xi_j - \eta}{\psi} \right).
\end{equation}
This density is known in the literature as the hyperbolic secant distribution. The horseshoe prior (\ref{eq:HS:xi}) is a special case of (\ref{eq:hyp:sech}) for $\eta=0$ and $\psi=1$. Without any loss of generality we can let $\eta=0$, as the location of the density will be determined by the value of the global shrinkage parameter $\tau$ as discussed in Section \ref{sec:log:scale:shrinkage:xi}. Allowing $\psi>0$ leads to a prior distribution for $\lambda_j = e^{\xi_j}$ of the form
\begin{equation}
	\label{eq:hyp_sech}
	p_{\rm SECH}(\lambda_j \vbar \psi) = \frac{ 2 \lambda_j^{1/\psi - 1}} { \pi \psi \left( \lambda_j^{2/\psi} + 1 \right) }.
\end{equation}
Examining (\ref{eq:hyp_sech}) clearly shows that the log-scale parameter $\psi$ controls the tail and concentration behaviour of the induced distribution over $\lambda_j$. The larger the log-scale $\psi$, the heavier the tail as $\lambda_j \to \infty$, and the greater the concentration of prior mass around $\lambda_j=0$. For $\psi=1$ this prior reduces to the half-Cauchy distribution; for $\psi<1$ the prior tends to zero as $\lambda_j \to 0$ and for $\psi>1$ the prior exhibits a pole at $\lambda_j=0$.

It is of interest to compare the prior (\ref{eq:hyp_sech}) to the prior over $\lambda_j$ one would obtain by starting with a $z$-distribution with shape parameters $a$ and $b$, mean of zero and a scale of $s>1$ over $\xi_j$ and transforming this to $\lambda_j$, which yields
\begin{equation}
	\label{eq:z:lambda:distr}
	p_{Z}(\lambda_j \vbar a,b,\psi) \propto \frac{\lambda_j^{2 a/s - 1}}{\left(1 + \lambda_j^{2/s} \right)^{a+b}}.
\end{equation}
Comparing (\ref{eq:z:lambda:distr}) with (\ref{eq:hyp_sech}) we see that the $a$ and $b$ shape parameters play exactly the same role as the log-scale parameter $\psi$, the primary difference being the ability to varying the tail or concentration behaviour individually by appropriate choice of $a$ and $b$. If we consider the case in which the shape parameters are the same, i.e., $a=b$, for which the $z$-distribution is symmetric, we have the following result.

\begin{proposition}
\label{prop:z:bound}
There exists a $K \in (0,\infty)$ such that
\[
	p_{\rm Z}(\lambda_j \vbar a,a,s) \leq K \, p_{\rm SECH}(\lambda_j \vbar \psi = s/(2 a))
\]
where $a>0$ is a shape parameter. 
\end{proposition}

The proof follows in a straightforward manner by application of the bound $1+x^c \leq (1+x)^c$. From the monotone convergence theorem, Proposition \ref{prop:z:bound} tells us that controlling the tails of the beta prime prior distribution over $\lambda_j$ by variation of the shape parameters cannot lead to a heavier tailed marginal distribution over $\beta_j$, or one with a greater concentration of mass at $\beta_j = 0$, than controlling the tails of the prior by varying the scale parameter $\psi$ of the log-$z$ distribution (\ref{eq:BP:tj}) alone. This result is also confirmed by the fact that the tails of $z$-distribution are log-linear with an absolute log-gradient of $a/s$ (\cite{BarndorffNielsenEtAl82}, p. 150). 


\subsection{Log-Laplace Priors}
\label{sec:LogLaplace}

We now examine a specific choice of log-scale prior based on the Laplace (double exponential) distribution. The primary usefulness of this distribution is its ability to provide simple bounds for the entire class of log-location-scale prior distributions over $\xi_j$ with log-linear, or sub-log-linear tails, which itself includes the important sub-class of log-concave densities. The log-Laplace prior distribution for a local shrinkage hyperparameter $\lambda_j$ is given by
\[
	\log \lambda_j \vbar \psi_1,\psi_2 \sim {\rm DE}(\psi_1,\psi_2),
\]
where ${\rm DE}(\psi_1,\psi_2)$ denotes an asymmetric Laplace distribution with a median of zero, a left-scale of $\psi_1$, a right-scale of $\psi_2$ and probability density function
\begin{equation}
	\label{eq:LogLaplace:Density:xi}
	p(\xi \vbar \psi_1,\psi_2) = \left( \frac{1}{2 \psi(\xi_j)} \right) \exp \left( -\frac{|\xi_j|}{\psi(\xi_j)} \right).
\end{equation}
where $\psi(\xi_j) = I(\xi_j<0)\psi_1 + I(\xi_j\geq 0)\psi_2$. The asymmetric Laplace distribution is essentially equivalent to two back-to-back exponential distributions with different scale parameters for each of the exponential distributions, and leads to a piecewise probability density function over $\lambda_j = e^{\xi_j}$ of the form:
%
%
%
\begin{equation}
	\label{eq:LogLaplacePDF}
	p_{\rm LL}(\lambda_j \vbar \psi_1, \psi_2) = \left\{\begin{array}{cc} 
			\displaystyle{ \left( \frac{1}{2 \psi_1} \right) \lambda_j^{\displaystyle{-1 + 1 / \psi_1}} } & \quad 0 < \lambda_j \leq 1 \\
		  \displaystyle{ \left( \frac{1}{2 \psi_2} \right) \lambda_j^{\displaystyle{-1-1/\psi_2}} } & \quad \lambda_j > 1 \\
		\end{array} \right..
\end{equation}
This distribution has a non-differentiable point at $\lambda_j=1$, and is discontinuous at $\lambda_j=1$ if $\psi_1 \neq \psi_2$. The piece of the function for $\lambda_j \in (0, 1)$ is proportional to a beta distribution, ${\rm Be}(\psi^{-1},1)$, and the piece of the function for $\lambda_j > 1$ is a Pareto distribution with shape parameter $1/\psi_2$. In the special case that $\psi_1 = \psi_2 = \psi$ the distribution reduces to the usual symmetric double exponential distribution which we denote by $p_{\rm LL}(\lambda_j \vbar \psi)$. An important property of the Laplace distribution is that it provides an upper-bound for the entire class of log-concave probability distributions. 

%
%

\begin{proposition}
\label{prop:log:concave:bounds}
Let $f(\xi_j)$ be a log-concave distribution with mode at $\xi_j = \xi^\prime$, and let $\xi_1<\xi^\prime$ and $\xi_2>\xi^\prime$ be any two values of $\xi$ on either side of $\xi^\prime$. Then, there exists a constant $K>0$ depending on $\xi^\prime$, $\xi_1$, $\xi_2$ such that
\[
	 f(\xi) \leq K \, p_{\rm LL}(\xi \vbar -g(\xi_1)^{-1}, g(\xi_2)^{-1}) \; \; \mbox{for all } \xi \in \mathbb{R},
\]
where $g(\xi) = -d \log f(\xi)/d \xi$ is the derivative of the negative logarithm of the density $f(\xi)$.
\end{proposition}

This is simply a restatement of a well known result regarding log-concave functions \citep{GilksWild92}. This result provides a useful upper bound which we use in Section \ref{sec:TheoreticalResults} to provide results regarding the concentration properties and tail behaviour of the entire class of log-concave prior distributions over $\xi_j$. For the specific case of the hyperbolic secant prior (\ref{eq:hyp_sech}) we can construct the following upper and lower bound based on the symmetric Laplace distribution.

\begin{proposition}
\label{prop:hyp:sec:bounds}
The log-hyperbolic secant distribution satisfies
\[
	\left( \frac{\pi}{2} \right) p_{\rm LL}(\xi_j \vbar \psi) \leq p_{\rm HS}(\xi_j \vbar \psi) \leq \left( \frac{4}{\pi} \right) p_{\rm LL}(\xi_j \vbar \psi)
\]
for all $\psi>0$.
\end{proposition}

The fact that the Laplace distribution on which these bounds are based has the same scale as the hyperbolic secant distribution it is bounding can be used to demonstrate that the log-Laplace distribution with scale $\psi$ and the log-HS distribution with scale $\psi$ lead to marginal distributions for $\beta_j$ that have identical concentration properties and tail behaviour. More generally, we have the following result.
\begin{proposition}
	\label{prop:log:linear:tails}
	Let $f(\xi)$ be a distribution over $\mathbb{R}$ that is bounded from above. If $f(\xi)$ satisfies
	\begin{equation}
		\label{eq:log:linear:tail:ass:1}
		f(\xi) = O \left( e^{a \xi} \right) \; \; \mbox{as } \xi \to -\infty, \; \; \mbox{and } f(\xi) = O \left( e^{-b \xi} \right) \; \; \mbox{as } \xi \to \infty
	\end{equation}
  where $a>0$, $b>0$, then there exists a $K_1>0$ such that
	\[
		f(\xi) \leq K_1 \, p_{\rm LL}(\xi \vbar a, b) \; \; \mbox{for all } \xi \in \mathbb{R}
	\]
	If $f(\xi)$ satisfies
	\begin{equation}
		\label{eq:log:linear:tail:ass:2}
		f(\xi) = \Omega \left( e^{a \xi} \right)  \; \; \mbox{as } \xi \to -\infty, \; \; \mbox{and } f(\xi) = \Omega \left( e^{-b \xi} \right) \; \; \mbox{as } \xi \to \infty
	\end{equation}
  where $a>0$, $b>0$, there exists a $K_2>0$ such that
	\[
		K_2 \, p_{\rm LL}(\xi \vbar a, b) \leq f(\xi)  \; \; \mbox{for all } \xi \in \mathbb{R}
	\]

\end{proposition}
This proposition tells us that the log-Laplace distribution can provide an upper bound for any distribution over $\xi$ which is log-linear, or sub-log-linear, in its tails, and can provide a lower-bound if the distribution is log-linear, or super-log-linear, in its tails. The advantage of these bounds are that the form of the log-Laplace distribution allows for relatively simple analysis of concentration and tail properties of the marginal distribution $p_{\rm LL}(\beta_j \vbar \psi_1,\psi_2)$, which we can use to derive bounds on the behaviour of the entire class of bounded prior densities over $\xi$ which have log-linear tails.

\subsection{The log-$t$ prior}
\label{sec:Logt}

The log-Laplace prior discussed in Section \ref{sec:LogLaplace} is important as it offers an upper-bound on the entire class of prior distributions on $\xi_j$ with log-linear tails, and through the monotone convergence theorem, an upper-bound on the marginal distributions over $\beta_j$ that they induce. It is of some interest then to examine an example of a prior distribution that cannot be bounded by the log-Laplace distribution. Specifically, we examine the log-$t$ prior distribution for $\xi_j = \log \lambda_j$:
\[
	\xi_j \sim t_\alpha(\psi)
\]
where $t_{\alpha}(\psi)$ denotes a Student-$t$ distribution centered at zero, with degrees-of-freedom $\alpha$, scale $\psi$ and probability density
\begin{equation}
	\label{eq:t:xi}
	p_t(\xi_j \vbar \alpha, \psi) = \frac{\Gamma((\alpha+1)/2)}{\Gamma(\alpha/2) \sqrt{\pi \alpha} \psi} \left( 1 + \frac{\xi_j^2}{\alpha \psi^2} \right)^{-(\alpha+1)/2}.
\end{equation}
Transforming the density (\ref{eq:t:xi}) to a density on $\lambda_j$ yields
\begin{equation}
	\label{eq:t:lambda}
	p_t(\lambda_j \vbar \alpha, \psi) = \left( \frac{\Gamma((\alpha+1)/2)}{\Gamma(\alpha/2) \sqrt{\pi \alpha} \psi} \right) \lambda_j^{-1} \left( \frac{\log(\lambda_j)^2}{\alpha \psi^2} + 1 \right)^{-(\alpha+1)/2}
\end{equation}
The density (\ref{eq:t:lambda}) is of the form $c \, \lambda_j^{-1} L_t(\lambda_j)$, where 
\begin{equation}
	\label{eq:L:t}
	L_t(x) = \left( \frac{\log(x)^2}{\alpha \psi^2} + 1 \right)^{-(\alpha+1)/2}
\end{equation}
is a function of slow variation~(\cite{BarndorffNielsenEtAl82}, p. 155). In this sense, the density (\ref{eq:t:lambda}) can be thought of as the normal-Jeffreys' prior $\lambda_j^{-1}$ multipled by a factor $L_t(\lambda_j)$ that slows its growth as $\lambda_j \to 0$, and increases the rate at which it decays as $\lambda_j \to \infty$, by an amount sufficient to ensure the resulting prior density is proper. The log-$t$ density dominates the log-Laplace density in the following sense.

\begin{proposition}
	\label{prop:t:dominate:LL}
	For all $\psi>0$, $\alpha>0$, $\psi_1>0$ and $\psi_2>0$ the log-$t$ density (\ref{eq:t:lambda}) satisfies 
	\[
		\lim_{\lambda_j \to 0^{+}} \left\{ \frac{p_t(\lambda_j \vbar \alpha, \psi)}{p_{\rm LL}(\lambda_j \vbar \psi_1,\psi_2)} \right\} = \infty \; \; \; \mbox{and} \; \; \;
		\lim_{\lambda_j \to \infty} \left\{ \frac{p_t(\lambda_j \vbar \alpha, \psi)}{p_{\rm LL}(\lambda_j \vbar \psi_1,\psi_2)} \right\} = \infty
	\]
	where $p_{\rm LL}(\cdot)$ is the log-Laplace density (\ref{eq:LogLaplacePDF}).
\end{proposition}

This result shows that the log-$t$ density, irrespective of the choice of degrees-of-freedom or scale parameter, always concentrates more probability mass near $\lambda_j=0$, and decays more slowly as $\lambda_j$ becomes large, than {\em any} prior density for $\lambda_j$ derived from a density on $\xi_j$ with log-linear tails. The log-$t$ density over $\xi_j$ implies a density over the shrinkage factor, $\kappa_j = 1/(1+\lambda_j^2)$ of the form
\[
	p_t(\kappa_j \vbar \alpha, \psi) \propto \kappa_j^{-1} (1-\kappa_j)^{-1} \left( 1 + \frac{{\rm arctanh}(1-2\kappa_j)^2}{\alpha \psi^2} \right)^{-(\alpha+1)/2}.
\]
	For all $\psi>0$ and degrees-of-freedom $\alpha>0$ this density satisfies
	\begin{equation}
		\label{prop:logt:kappa}
		\lim_{\kappa_j \to 0^{+}} \left\{ p_t(\kappa_j \vbar \alpha, \psi) \right\} = \infty \; \; \; \mbox{and} \; \; \; \lim_{\kappa_j \to 1^{-}} \left\{ p_t(\kappa_j \vbar \alpha, \psi) \right\} = \infty,
	\end{equation}
	%
%
which shows that regardless of the choice of degrees-of-freedom parameter $\alpha$, or the scale parameter $\psi$, the log-$t$ prior distribution leads to a prior distribution over $\kappa_j$ that is always infinite at ``no shrinkage'' ($\kappa_j=0$) and ``complete shrinkage to zero'' ($\kappa_j=1$). However, despite this property, the log-$t$ density can concentrate as much probability mass around $\kappa_j=1/2$ as desired by an appropriate choice of $\psi$, as formalised by the following result.

\begin{proposition}
	\label{prop:logt:kappa:concentration}
	For all degrees-of-freedom $\alpha>0$, $\epsilon \in (0,1/2)$ and $\delta \in (0,1)$ there exists a $\psi>0$ such that
	\[
		\int_{1/2-\epsilon}^{1/2+\epsilon} p_t(\kappa_j \vbar \alpha, \psi) d \kappa_j > \delta
	\]
\end{proposition}

This result implies that by choosing a small enough log-scale $\psi$, the log-$t$ prior distribution can become more and more similar to the ridge regression prior by allowing less {\em a priori} variation in shrinkage between model coefficients. However, the property (\ref{prop:logt:kappa}) guarantees that regardless of how much prior probability mass is concentrated around $\kappa=1/2$ the density $p_t(\kappa_j \vbar \cdot)$ always tends to infinity as $\kappa_j \to 0$ and $\kappa_j \to 1$. 


\section{Discussion and Theoretical Results}
\label{sec:TheoreticalResults}

In this section we examine the theoretical behaviour of the log-scale prior distributions for $\lambda_j$, when used within the hierarchy
\[
	\beta_j \vbar \lambda_j \sim N(0, \lambda_j^2), \; \; \; \log \lambda_j \vbar \psi \sim f(\log \lambda_j) d \log \lambda_j,
\]
where $f(\cdot)$ is a unimodal distribution over $\log \lambda_j$. Define the marginal distribution $\pi_f(\beta_j)$ of $\beta_j$, relative to the prior distribution $f(\cdot)$, by
\begin{equation}
	\label{eq:MarginalForBeta}
	\pi_{f}(\beta_j) = \int_{0}^\infty \left( \frac{1}{2 \pi \lambda_j^2} \right)^{\frac{1}{2}} \exp \left( - \frac{\beta_j^2}{2 \lambda_j^2} \right) f(\lambda_j) d \lambda_j.
\end{equation}
As discussed in Section \ref{sec:intro:properties}, two desirable properties of a prior distribution $f(\lambda_j)$ over $\lambda_j$ are that corresponding marginal distribution $\pi_f(\beta)$: (I) tends to infinity as $|\beta_j| \to 0$, and (II) has Cauchy or super-Cauchy tails as $|\beta_j| \to \infty$. We will now show that for an appropriate choice of scale parameters the asymmetric log-Laplace prior distribution results in a marginal distribution that posesses both of these properties. First, we examine the form of the marginal distribution when $f(\lambda_j)$ is an asymmetric log-Laplace distribution.

\begin{theorem}
\label{thm:LLMarginal}
Let $\lambda_j$ follow a log-Laplace distribution with left log-scale $\psi_1$ and right log-scale $\psi_2$, and let $\beta_j$ follow a normal distribution with variance $\lambda_j^2$. The marginal distribution (\ref{eq:MarginalForBeta}) for the regression coefficient $\beta_j$ is
\begin{equation}
	\label{eq:MarginalForBeta:LL}
	\pi_{\rm LL}(\beta_j \vbar \psi_1,\psi_2) = \left( \frac{1}{\sqrt{32 \pi}} \right) \left[ \frac{1}{\psi_1} E_{\left(\frac{1+\psi_1}{2 \psi_1}\right)} \left( \frac{\beta_j^2}{2} \right) + \frac{1}{\psi_2} \left( \frac{2}{\beta_j^2} \right)^{\left(\frac{1+\psi_2}{2 \psi_2}\right)} \gamma \left( \frac{1+\psi_2}{2 \psi_2}, \frac{\beta_j^2}{2} \right) \right]
\end{equation}
where $E_n(\cdot)$ is generalized exponential integral and $\gamma(s,x)$ is the incomplete lower-gamma function.
\end{theorem}

Using Theorem \ref{thm:LLMarginal} we can examine the concentration properties of the marginal distribution when $f(\lambda_j)$ is an asymmetric log-Laplace distribution.

\begin{theorem}
\label{eq:LogLaplaceTheorem1}
Let $\lambda_j$ follow a log-Laplace distribution with left log-scale $\psi_1$ and right log-scale $\psi_2$, and let $\beta_j$ follow a normal distribution with variance $\lambda_j^2$. Then, for all $\psi_2>0$, as $|\beta_j| \to 0$, the marginal density satisfies
\begin{enumerate}
	\item $\pi_{\rm LL}(\beta_j \vbar \psi_1,\psi_2) = O \left( |\beta|^{-1+1/\psi_1} \right)$ if $\psi_1 > 1$;
	
	\item $\pi_{\rm LL}(\beta_j \vbar \psi_1,\psi_2) = O \left( -\log |\beta|_j \right)$ if $\psi_1 = 1$;
	
	\item $\pi_{\rm LL}(\beta_j \vbar \psi_1,\psi_2) = O(1)$ if $\psi_1 < 1$.
\end{enumerate}
as $|\beta_j| \to 0$.

\end{theorem}

We also have the following theorem, which characterises the tail behaviour of marginal distribution when $f(\lambda_j)$ is an asymmetric log-Laplace distribution. 

\begin{theorem}
\label{eq:LogLaplaceTheorem2}
Let $\lambda_j$ follow a log-Laplace distribution with left log-scale $\psi_1$ and right log-scale $\psi_2$, and let $\beta_j$ follow a normal distribution with variance $\lambda_j^2$. Then, for all $\psi_1>0$ 
\[
	\pi_{\rm LL}(\beta_j \vbar \psi_1,\psi_2) = O \left( |\beta|^{-1-1/\psi_2} \right)
\]
as $|\beta| \to \infty$
\end{theorem}

\begin{figure*}[t]
\begin{center}
\subfigure[$\psi=1/2$]{
   \includegraphics[scale=0.235]{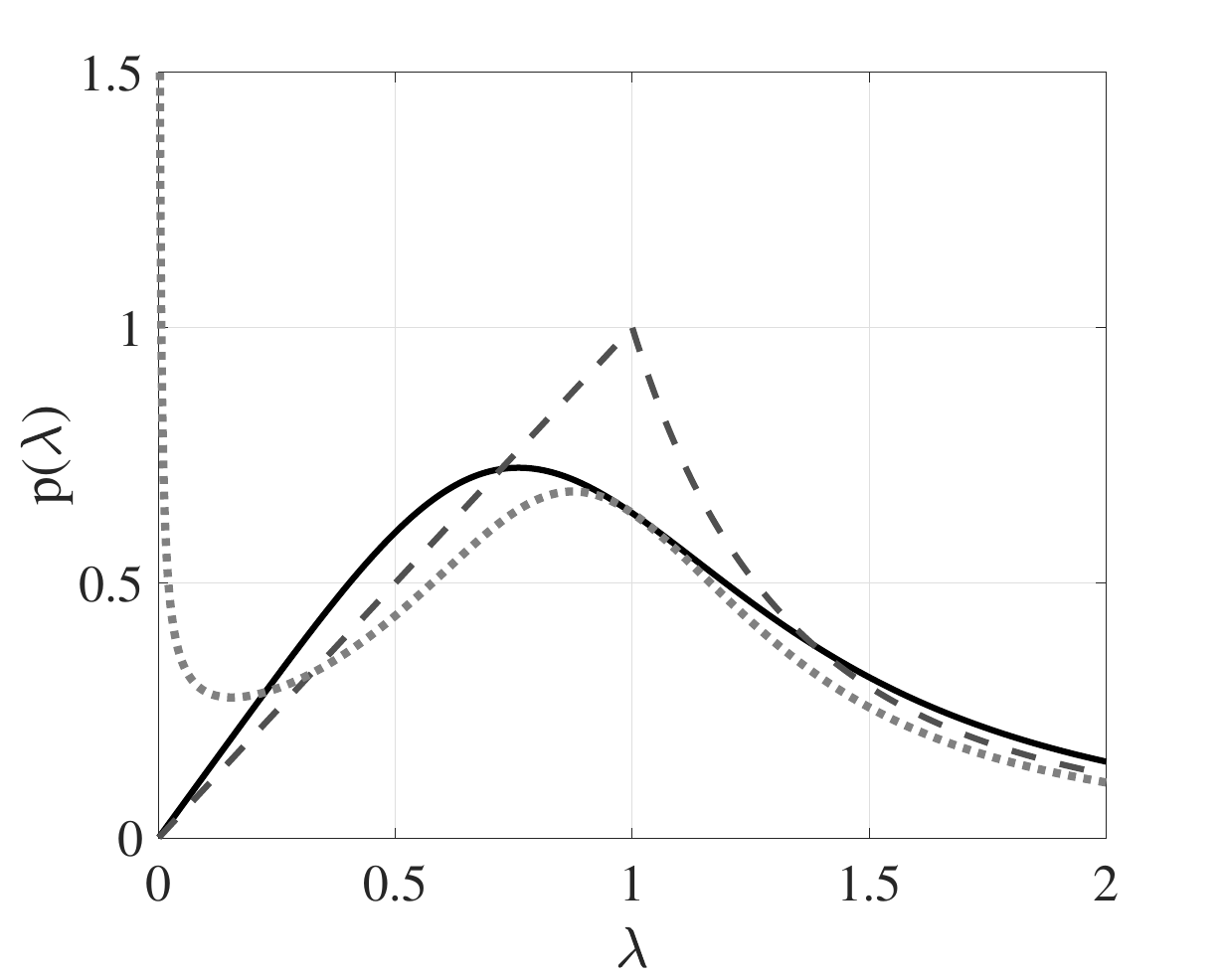}
}%
\subfigure[$\psi=1$]{
   \includegraphics[scale=0.235]{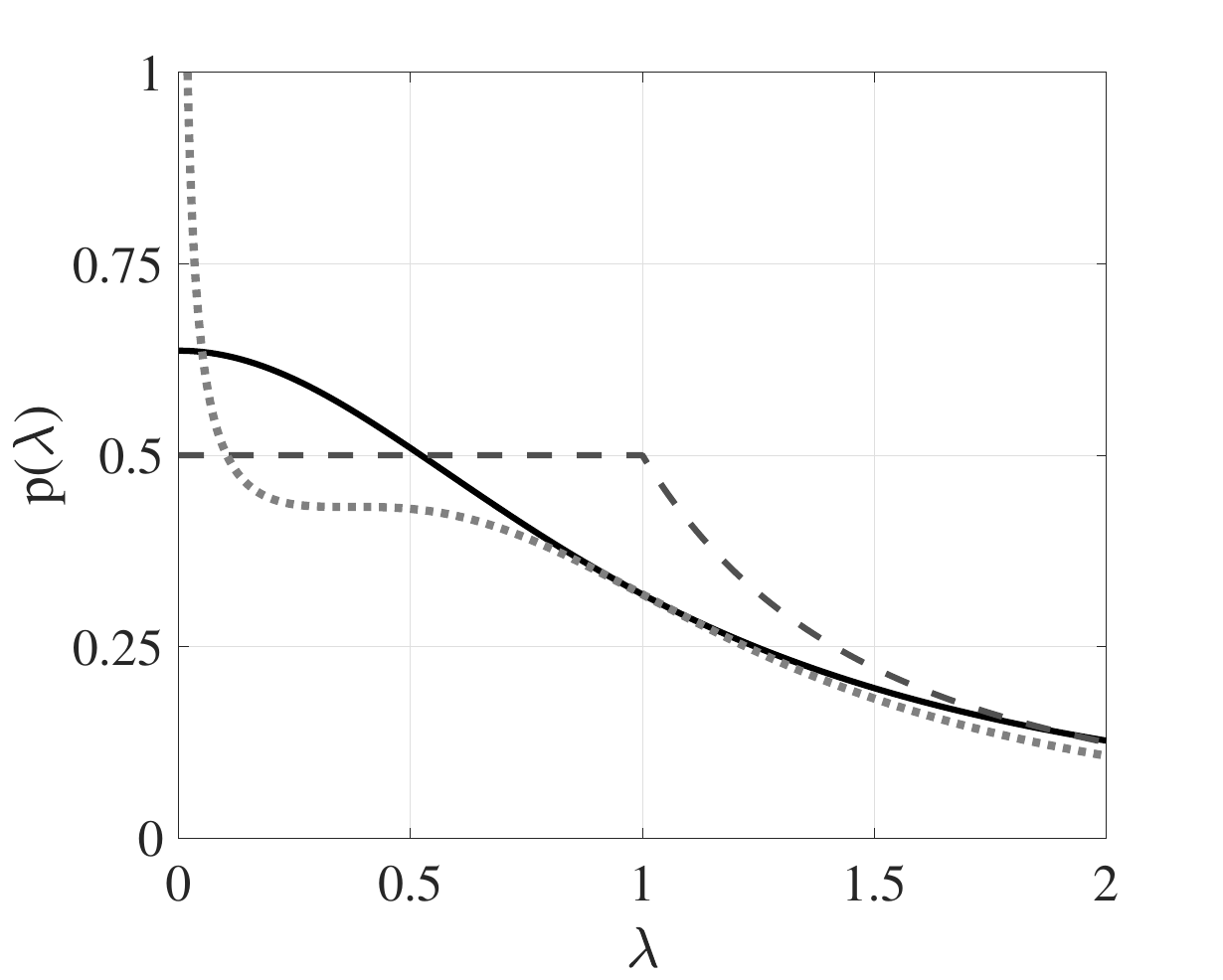}
}
\subfigure[$\psi=2$]{
   \includegraphics[scale=0.235]{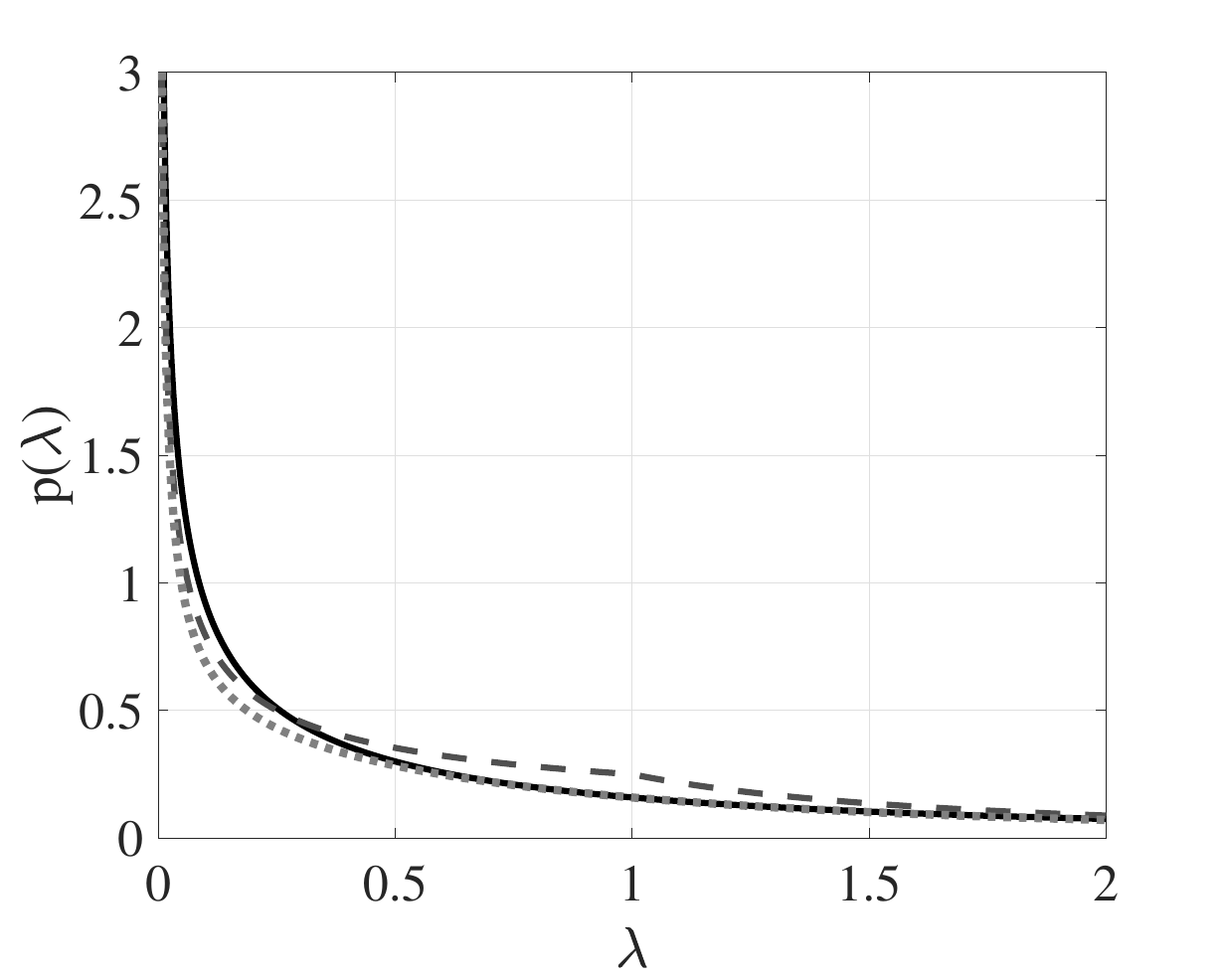}
}%
\\
\subfigure[$\psi=1/2$]{
   \includegraphics[scale=0.235]{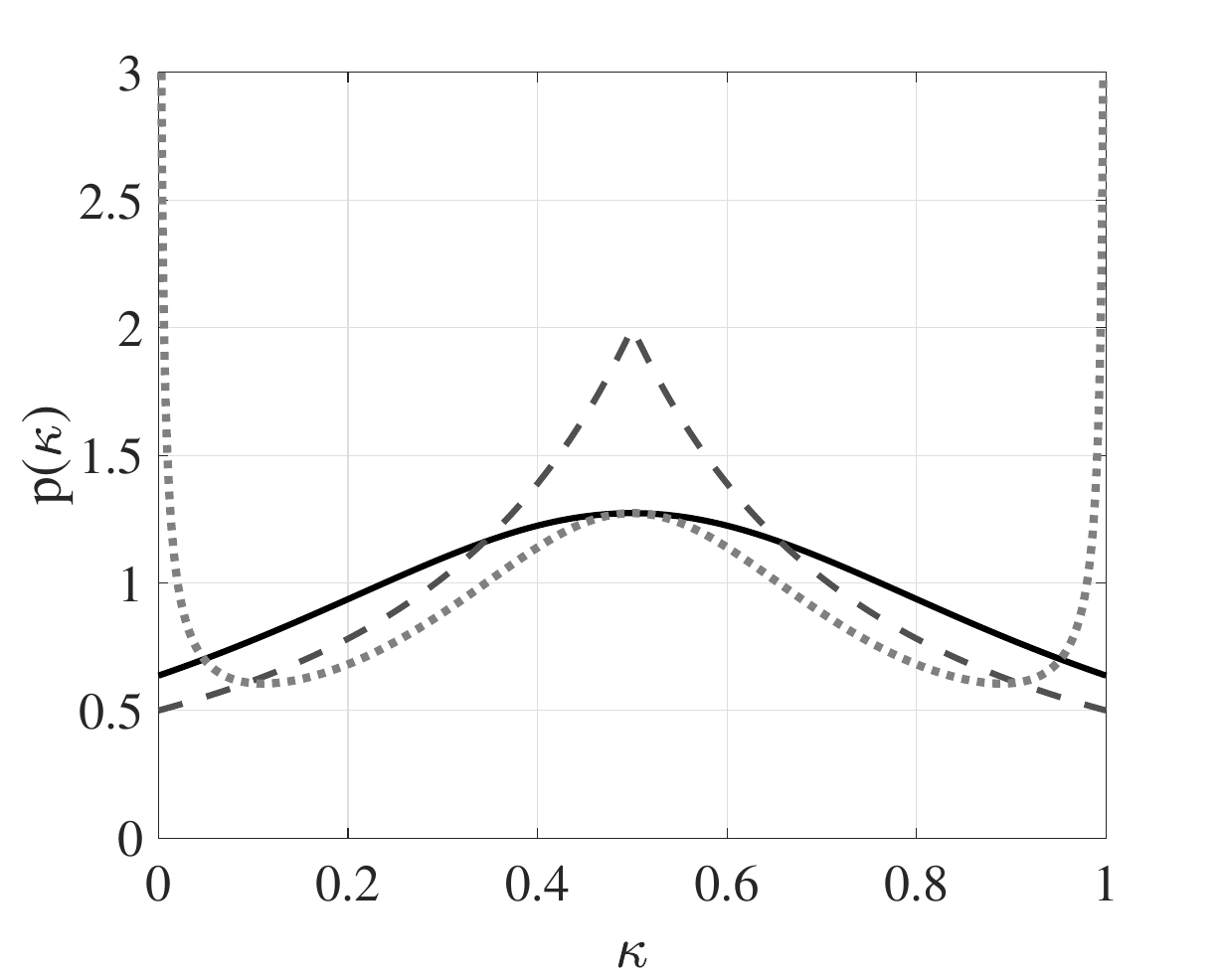}
}%
\subfigure[$\psi=1$]{
   \includegraphics[scale=0.235]{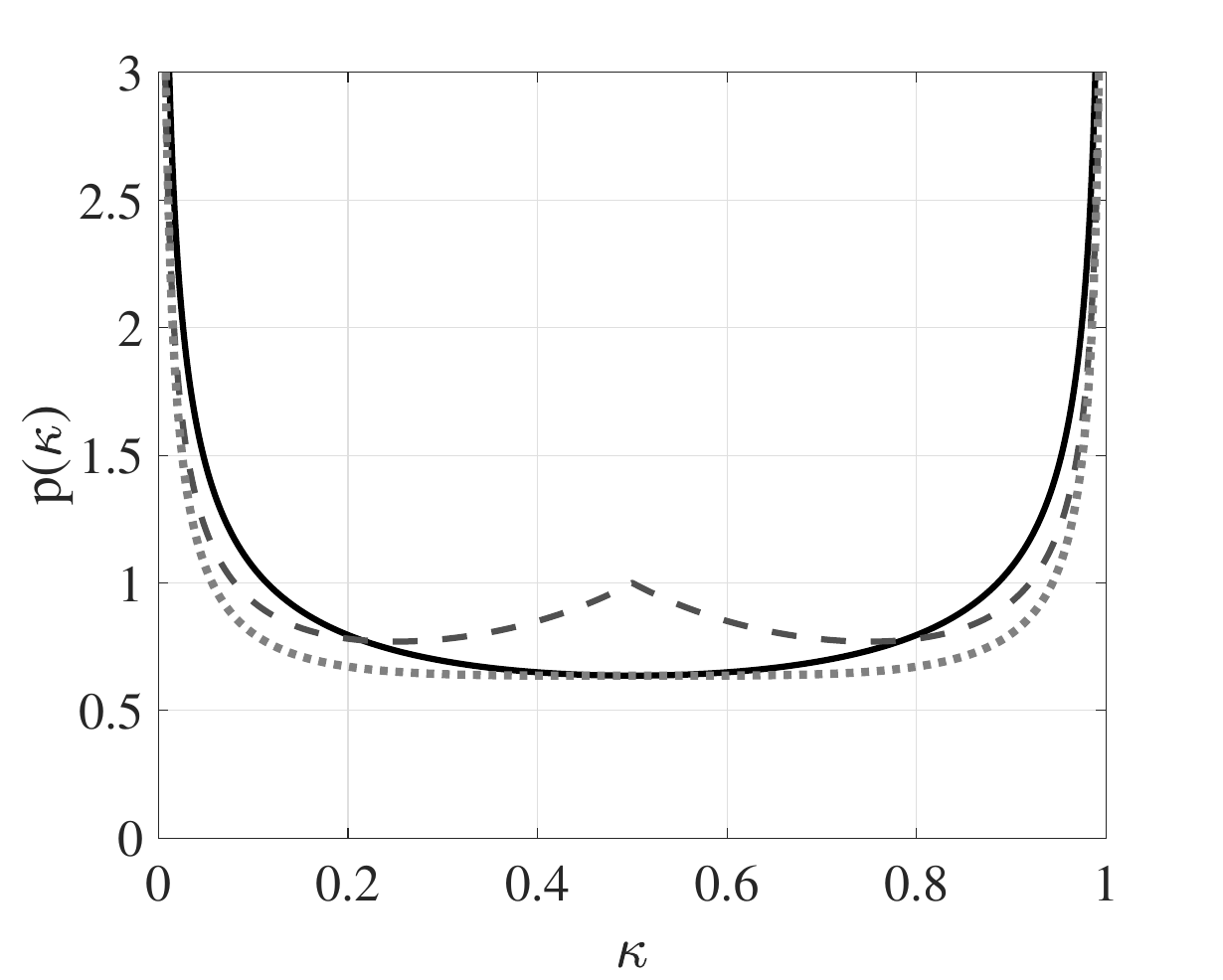}
}
\subfigure[$\psi=2$]{
   \includegraphics[scale=0.235]{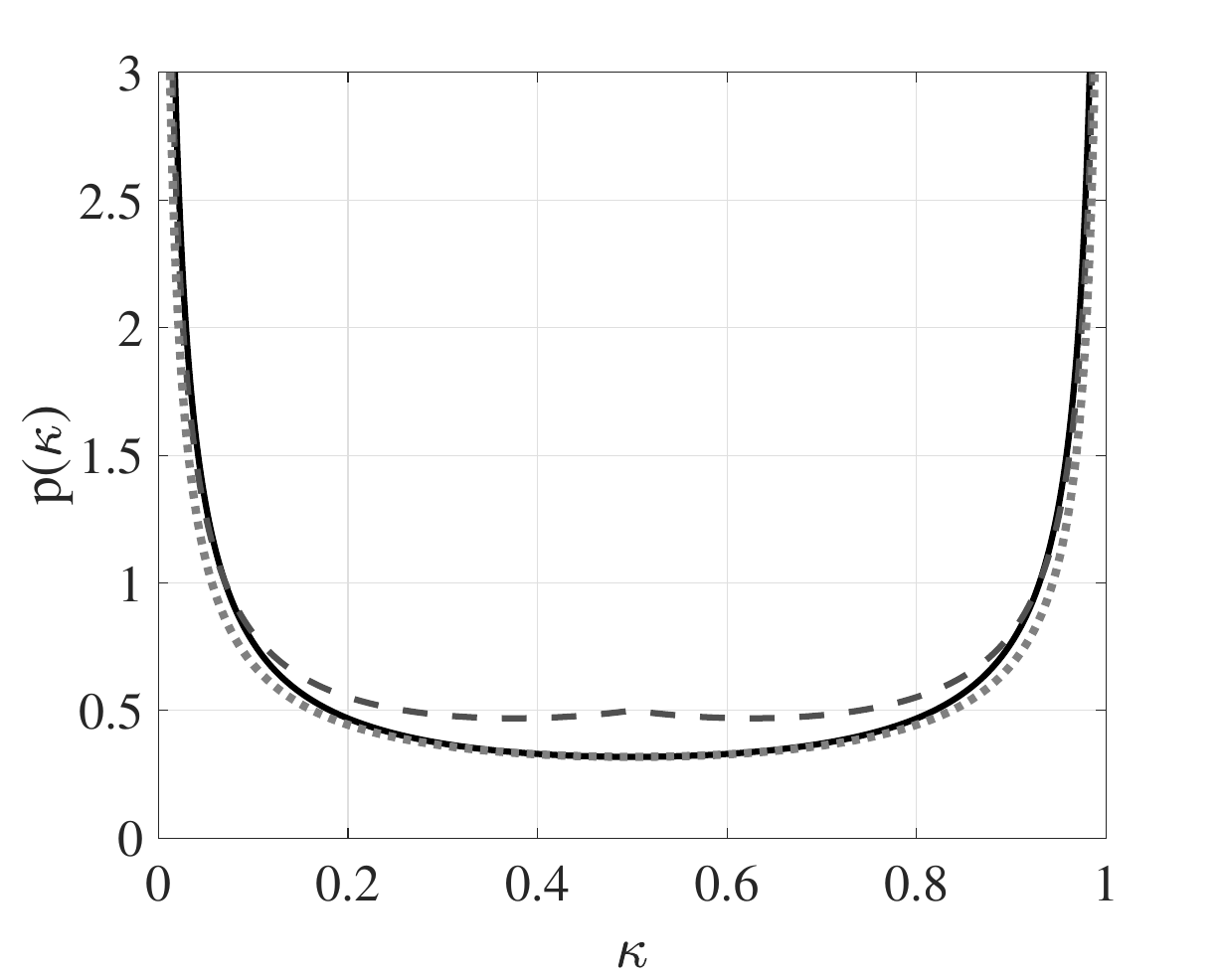}
}%
\end{center}
\caption{Prior probability density plots for the log-hyperbolic secant (solid), log-Laplace (dashed) and log-$t$ with $\alpha=1$ (dotted) distributions for the local shrinkage parameter $\lambda$, for three different choices of log-scale $\psi$, and the corresponding implied prior distributions on the coefficient of shrinkage $\kappa$. Note that when $\psi=1$, the log-hyperbolic secant distribution is equivalent to the usual horseshoe (half-Cauchy) prior distribution.\label{fig:priors:lambda}}
\end{figure*}

Theorems \ref{eq:LogLaplaceTheorem1} and \ref{eq:LogLaplaceTheorem2}, combined with Proposition \ref{prop:hyp:sec:bounds} and the monotone convergence theorem can be used to obtain identical results for the log-hyperbolic secant hyperprior. More generally, when combined with Proposition \ref{prop:log:linear:tails}, Theorems \ref{eq:LogLaplaceTheorem1} and \ref{eq:LogLaplaceTheorem2} provide upper bounds on prior concentration at $\beta_j=0$ for the entire class of prior distributions over $\xi_j$ with log-linear, or sub-log-linear tails. An important aspect of these results is that both priors exhibit a pole at $|\beta_j| = 0$, and have Cauchy or super-Cauchy tails as $|\beta_j| \to \infty$, if and only if $\psi_1, \psi_2 \geq 1$. If $\psi_1 < 1$, insufficient mass is concentrated at $\lambda_j$ to produce a pole at $\beta_j=0$. If $\psi_2 < 1$ the tail of the prior distribution over $\lambda_j$ is too light, and the marginal distribution of decays at a super-Cauchy rate. In contrast, when $\xi_j$ follows a Student-$t$ distribution, as per Section (\ref{sec:Logt}), the log-Laplace prior density no longer provides an upper-bound, and the resulting marginal distribution for $\beta_j$ exhibits very different behaviour, as characterised by the following result.


\begin{theorem}
\label{eq:Log:student:t:theorem}
Let $\lambda_j$ follow a log-$t$ distribution with log-scale $\psi$ and degrees-of-freedom parameter $\alpha$, and let $\beta_j$ follow a normal distribution with variance $\lambda_j^2$. Then, the resulting marginal distribution over $\beta_j$ satisfies 
\[
	\pi_{t}(\beta_j \vbar \alpha, \psi) = \Omega \left( |\beta_j|^{-1+1/c} \right) \; \; \mbox{as } |\beta_j| \to 0
\]
%
for all $c>0$, and
\[
	\pi_{t}(\beta_j \vbar \alpha, \psi) \asymp  |\beta_j|^{-1} (\log |\beta_j|)^{-\alpha-1} \; \; \mbox{as } |\beta_j| \to \infty
\]
%
for all $\psi>0$ and $\alpha>0$.
\end{theorem}

This result demonstrates a very interesting property of the log-$t$ prior distribution: namely, that irrespective of the choice of the degrees-of-freedom parameter $\alpha$, or the scale parameter $\psi$, the resulting marginal distribution over $\beta_j$ always possesses Properties I and II. To the best of the authors' knowledge, this property appears to be unique amongst all the known prior distributions for $\lambda_j$.

Figure \ref{fig:priors:lambda} shows plots of the log-hyperbolic secant, log-Laplace and log-$t$ (with degrees-of-freedom $\alpha=1$) prior distributions over the local shrinkage parameter $\lambda$. The manner in which the log-Laplace prior distribution serves as a general analog for prior distributions over $\xi = \log \lambda$ with log-linear tails is clear by visual comparison with the log-hyperbolic secant (in this case, the usual horseshoe prior). When $\psi = 1/2$, it is clear that both the log-hyperbolic secant and log-Laplace distributions tend to zero at $\lambda=0$ which results in the prior density for $\kappa$ tending to a constant as $\kappa \to 0$. In contrast, the log-$t$ prior distribution with $\psi = 1/2$ still tends to infinity as $\kappa \to 0$, as predicted by Proposition~\ref{prop:t:dominate:LL} (and will do so for any $\psi>0$), and by Theorem~\ref{eq:Log:student:t:theorem}, the resulting marginal density over $\beta$ will tend to infinity as $|\beta| \to 0$, and possess heavier tails than any member of the log-Laplace family (and by extension, any prior density over $\log \lambda$ with log-linear tails, such as the horseshoe and horseshoe+).

\subsection{Kullback--Leibler Super-efficiency}

We now examine the behaviour of the estimators implied by log-scale prior distributions in terms of their average accumulated Kullback--Leibler (Ces\`{a}ro) risk. let $\hat{p}_n(x)	 = \int p(x|\bm{\theta}) p(\bm{\theta}|y_1,\ldots,y_n) d \theta$ denote the posterior predictive density, and let $D(p || q) = \int p(x) \log p(x)/q(x) dx$ denote the Kullback--Leibler divergence from $p$ to $q$; the Ces\`{a}ro risk at $\bm{\theta}_0$ is then given by
\[
	R^\pi_n(\bm{\theta}_0) = \frac{1}{m} \sum_{n=1}^N \E{\bm{\theta}_0}{ D(p_{\bm{\theta}_0} \vbar \hat{p}_n) ) }
\]
where $\E{\beta}{\cdot}$ denotes expectation with respect to the sampling distribution $p(y_1,\ldots,y_n \vbar \bm{\theta}_0)$. The Ces\`{a}ro risk can be interpreted as the expected average accumulated prediction error obtained by using the information contained in the posterior distribution $\hat{p}_n$ based on $n$ samples to predict the $(n+1)$-th sample, for $n=1$ to $N$. The first tool that will be used is the following result from \cite{ClarkeBarron90}:

\begin{lemma}
\citep{ClarkeBarron90}. Let $p(A) = \int_A \pi(\bm{\theta}) d\bm{\theta}$ denote the prior mass assigned to set $A \subseteq \Theta$ by the prior distribution $\pi(\bm{\theta})$, and let $\pi(\bm{\theta})$ be information rich at every point $\bm{\theta}_0 \in \Theta$, in the sense that $p(A_{\bm{\theta}_0, \varepsilon}) > 0$, where $A_{{\bm{\theta}_0,\varepsilon}} = \{ \bm{\theta} : D(p_{\bm{\theta}_0} || \bm{\theta}) \leq \epsilon \}$ is an $\varepsilon$-radius KL ball centred at $\bm{\theta}_0$. Then, the Cesa\`{r}o-risk for the Bayes procedure based on $\pi(\bm{\theta})$ satisfies
\[
	R^\pi_n(\bm{\theta}_0) \leq \varepsilon -  \frac{1}{n} \log p(A_{{\bm{\theta}_0,\varepsilon}})
\]
for all $n>0$ and $\varepsilon > 0$.
\end{lemma}

This remarkable theorem allows us to characterise the behaviour of a Bayes procedure based on a prior distribution $\pi(\bm{\theta})$ as the sample size $n$ grows. In the context of the multiple means model (\ref{eq:multiple:means}), we modify the setting slightly so that we now observe a matrix of data ${\bf Y} = ({\bf y}_1, \ldots, {\bf y}_q)$ where column ${\bf y}_j = (y_{1,j},\ldots,y_{n,j})$ corresponds to the datapoints that correspond to mean $\beta_j$; in this formulation the multiple-means problem is essentially equivalent to orthogonal regression. As the KL risk for the multiple-means problem is additive in the co-ordinates, and each of the co-ordinates is exchangable, we concentrate on the Cser\`{a}ro risk for a single co-ordinate. In this case, \cite{CarvalhoPolson10} showed that the risk upper bound for the least-squares is given by $R_n(\beta_0) = O \left(n^{-1} \log n \right)$ for all $\beta_0 \in \mathbb{R}$.

To demonstrate that the log-scale prior distributions are Kullback--Leibler super-efficient, in the sense that $R_n(\beta_0)$ tends to zero at a quicker rate than the least squares estimator when $\beta_j = 0$ we once again utilise the properties of the log-Laplace distribution; specifically, we use the ability of log-Laplace distribution to provide bounds for a large class of functions. To this end, we introduce the double log-Laplace distribution with parameters $\psi_1, \psi_2>0$, which has the probability density function
\[
	p_{\rm DLL}(\beta_j \vbar \psi_1,\psi_2) = \left( \frac{1}{2} \right) p_{\rm LL}(|\beta_j| \vbar \psi_1, \psi_2),
\]
with $p_{\rm LL}(\cdot)$ given by (\ref{eq:LogLaplacePDF}). The double log-Laplace distribution is formed by placing two log-laplace distributions back-to-back, and is symmetric around $\beta_j = 0$. The usefulness of the double log-Laplace is that it may act as an appropriate lower and upper bound for a wide range of prior distributions over $\beta_j$. Specifically, we can establish the following Ces\`{a}ro risk bounds for prior distributions that can be lower-bounded by an appropriate double log-Laplace distribution. 

\begin{theorem}
\label{thm:KL:super:eff}
Let $\pi(\beta)$ denote a prior distribution over $\beta \in \mathbb{R}$ that satisfies the following conditions: 
\begin{enumerate}
	\item $\pi(\beta) = \Omega \left( |\beta|^{-1+1/\psi_1} \right)$ as $|\beta| \to 0$ where $\psi_1 > 0$;
	\item $\pi(\beta) = \Omega \left( |\beta|^{-1-1/\psi_2} \right)$ as $|\beta| \to \infty$ where $\psi_2 > 0$; and
	\item $\pi(\beta) \in (0, \infty)$ for $|\beta| \in (0, \infty)$. 
\end{enumerate}
Then, if $n \psi_1 \geq 1$, the Ces\`{a}ro risk at $\beta_0 = 0$ for the Bayes procedure based on the prior $\pi(\cdot)$ and sampling distribution $y_1,\ldots,y_n \sim N(\beta, 1)$ satisfies 
\begin{eqnarray*}
	R_n^\pi(0) &\leq& \frac{1}{n} \left[ \left( \frac{1}{2 \psi_1} \right) \left( \log ( n \psi_1 ) + 1 \right) - \log K/2 \right] \\
	&=& \frac{1}{n} \left( 1 + \frac{\log n}{2} - \left( \frac{1}{2} - \frac{1}{2 \psi_1} \right) \log n + O(1) \right)
\end{eqnarray*}
where $K > 0$ is a constant independent of the sample size $n$.
\end{theorem}

By using Theorem \ref{thm:KL:super:eff} in conjunction with Theorem \ref{eq:LogLaplaceTheorem1} we can easily establish that the log-Laplace prior distribution (\ref{eq:LogLaplacePDF}) over the local shrinkage parameter $\lambda$ results in Kullback--Leibler superefficiency at $\beta_0 = 0$ as long as $\psi_1 > 1$. Additionally, using Proposition \ref{prop:log:concave:bounds} allows us to extend this result to the class of prior distributions on $\lambda$ that induce a prior distribution on $\xi = \log \lambda$ that have tails that decay at a super-linear rate, which includes the log-hyperbolic secant prior discussed in Section \ref{sec:LogHyperbolicSecant}; as discussed in the next section, the log-Laplace prior distribution lower-bounds the majority of the standard local shrinkage parameter prior distributions, which makes Theorem \ref{thm:KL:super:eff} a very general tool. Finally, from Proposition \ref{prop:t:dominate:LL}, it is clear that the log-$t$ prior distribution for $\lambda$ results in Kullback--Leibler super-efficiency at $\beta_0 = 0$ for all choices of the hyperparameters $\alpha>0$, $\psi>0$.

%
%

\subsection{Comparison with Standard Shrinkage Priors}
\label{sec:Comparison:Standard:Priors}

It is interesting to compare the new log-Laplace and log-$t$ prior distributions discussed in Sections \ref{sec:Logt} and \ref{sec:LogLaplace} with the standard shrinkage priors proposed in the literature. As a consequence of Proposition \ref{prop:log:concave:bounds}, the log-Laplace density is useful because it serves as an upper-bound for the entire class of probability densities that are log-concave on $\xi$; therefore, no prior density over $\xi$ that is log-concave can achieve greater concentration of marginal prior probability around $\beta=0$, or heavier tails as $|\beta| \to \infty$, than some member of the log-Laplace family. It is therefore interesting to determine which of the standard shrinkage prior distributions from the literature fall into this class.

The Bayesian lasso prior (\ref{eq:BayesLassoXi}) and the regular horseshoe prior (\ref{eq:HS:xi}) over $\xi$ are easily verified to be log-concave by examination of their second derivatives. More generally, the beta prime family of prior densities over $\lambda^2$ are characterised by the $z$-distribution (\ref{eq:BP:tj}) over $\xi$, which is log-concave. This implies that regardless of how the shape hyperparameters $a$ and $b$ are chosen, the beta prime prior cannot result in a marginal distribution for $\beta$ with greater concentration near $\beta=0$, or heavier tails as $|\beta|\to\infty$, than a log-Laplace density (\ref{eq:LogLaplacePDF}) with appropriately chosen scale parameters.

A recent trend has been to propose prior densities for $\lambda$ of the form
\[
	\lambda = \phi_1 \phi_2, \; \; \; \phi_1 \sim p(\phi_1)d\phi_1, \; \; \; \phi_2 \sim p(\phi_2)d \phi_2.
\]
The horseshoe+, R2-D2 prior~\citep{ZhangEtAl16} and inverse gamma-gamma (IGG) prior~\citep{BaiGhosh17a} can all be represented in this manner. To prove that the density induced by this prior over $\xi = \log \lambda = \log \phi_1 + \log \phi_2$ is log-concave, it suffices to show that $\log \phi_1$ and $\log \phi_2$ are both distributed as per log-concave densities; the preservation of log-concavity under convolution therefore guarantees that the distribution of their sum will also be log-concave (\cite{SaumardWellner14a}, pp. 60--61). In the case of horseshoe+, $\log \phi_1$ and $\log \phi_2$ are both distributed as per hyperbolic secant distributions (\ref{eq:HS:xi}), which are log-concave. The IGG prior is given by
\[
	\phi_1 \sim {\rm IG}(a,c), \; \; \; \phi_2 \sim {\rm Ga}(b, c).
\]
Irrespective of the choice of hyperparameters $a$, $b$ and $c$ both $\log \phi_1$ and $\log \phi_2$ are distributed as per log-concave densities, and it follows again that the IGG shrinkage prior for $\xi$ is log-concave, and its behaviour is also bounded by the log-Laplace shrinkage prior. Finally, the R2-D2 prior (see Equation 7, \cite{ZhangEtAl16}) is built as a scale mixture of double-exponential prior distributions over $\beta$:
\[
	\beta \vbar \lambda \sim {\rm La}(\sqrt{2} \lambda), \; \; \; \lambda^2 \sim {\rm BP}(a, b)
\]
${\rm BP}(\cdot)$ is the beta-prime density, and $a>0$ and $b>0$ are hyperparameters that control the tails of the density over $\lambda_j$. Using the standard scale-mixture of normals representation of the double-exponential distribution and we can rewrite this hierarchy as
\[
	\beta \vbar \phi_1, \phi_2 \sim N(0, \phi_1^2 \phi_2^2), \; \; \; \phi_1^2 \sim {\rm Exp}(1), \; \; \; \log \phi_2 \sim Z(0,a,b,1/2)
\]
from which it immediately follows that both $\log \phi_1$ and $\log \phi_2$ are distributed as per log-concave probability distributions. The interesting result here is that choice to to use a double-exponential kernel for $\beta$ in place of a normal kernel made in \cite{ZhangEtAl16} was motivated by the aim of producing a marginal prior density for $\beta$ with greater concentration at the origin, and heavier tails. However, as the implied density over $\xi$ is log-concave it is clear that the use of the double-exponential kernel does not lead to prior with any different asymptotic properties than simply using a normal kernel with an appropriate log-Laplace prior over $\lambda$.

In contrast to the above prior distributions, the log-$t$ distribution does not have (sub) log-linear tails, and therefore does not fit into the class of prior distributions upper-bounded by the log-Laplace prior. However, our work is not the first to propose a larger, unifying class of shrinkage priors. In particular, recent work by~\cite{GhoshChakrabarti17} has explored the consistency properties of a class of global-local shrinkage priors that can be written in a decomposition of the form
\begin{equation}
	\label{eq:ghosh:chakra:p}
 p(\lambda_j^2) = c \, \left( \lambda_j^2 \right)^{-a-1} L(\lambda_j^2)
\end{equation}
where $a>0$ and $L(\lambda_j^2)$ is a slowly-varying function that satisfies $\lim_{t \to \infty} L(t) \in (0,\infty)$. To show that the log-$t$ prior distribution also falls outside of this class of priors, we can transforming the log-$t$ distribution (\ref{eq:t:lambda}) over $\lambda_j$ to a distribution over $\lambda_j^2$, yielding
\begin{equation}
	\label{eq:p:t:lambda2}
	p_t(\lambda_j^2) = c \, \lambda_j^{-2} L_t(\lambda_j^2)
\end{equation}
where $L_t(\cdot)$ is given by (\ref{eq:L:t}). To express (\ref{eq:p:t:lambda2}) in the form (\ref{eq:ghosh:chakra:p}) we require that $a=0$. Further, we note that $L_t(x)$ tends to zero as $x \to \infty$. The log-$t$ prior therefore violates both of the conditions required for a prior distribution to fall into the particular class studied by~\cite{GhoshChakrabarti17}.

\subsection{Concentration Inequalities}

We now develop several inequalities regarding the concentration of posterior probability in subsets of $\kappa$-space. Similar inequalities have been derived by~\cite{DattaGhosh13} and ~\cite{GhoshEtAl16}. The developments are different and express the inequalities directly in terms of the cumulative frequency distribution of the underlying log-scale prior distribution, which is usually known. We examine the distribution of a single $\kappa_j = 1/(1+\lambda_j)$ hyperparameter conditional on both $y_j$ and any relevant hyperparameters such as $\tau$ and $\psi$. We utilise the fact that the standard global-local shrinkage hierarchy (\ref{eq:Global}) can be rewritten as
\begin{eqnarray}
	\nonumber
	y_j \vbar \beta_j &\sim& N(\beta_j, 1) \\
	\label{eq:GL:Hierarchy:2}
	\beta_j \vbar \lambda_j &\sim& N(0, e^{2 \xi_j}) \\
	\nonumber
	\xi_j \vbar \tau, \psi &\sim& p(\xi_j - \log \tau \vbar \psi)
\end{eqnarray}
in which the global shrinkage hyperparameter is moved up the hierarchy to act as a location parameter for $\xi_j = \log \lambda_j$. Following~\cite{BhadraEtAl15}, we marginalise out $\beta_j$, and write the posterior distribution of $\kappa_j = 1/(1+\lambda_j^2)$ as
\begin{eqnarray}
	\label{eq:posterior:kappa}
	p(\kappa_j \vbar y_j, \tau, \psi) = c(y_j, \tau, \psi) \kappa_j^{1/2} e^{-\kappa_j y_j^2/2} p(\kappa_j \vbar \tau, \psi)
\end{eqnarray}
where $p(\kappa_j \vbar \tau, \psi)$ is the implied prior distribution over the coefficient of shrinkage. Writing $p(\xi_j - \log \tau \vbar \psi) \equiv p(\xi_j \vbar \log \tau, \psi)$, the following proposition gives general bounds on the posterior probability that $\kappa$ lies on edges of the $\kappa$-space; that is, the posterior probability that the shrinkage coefficient is close to zero or one.

\begin{proposition}
\label{prop:kappa:ineq}
Let $p(\xi_j \vbar \log \tau, \psi)$ be a unimodal, fully-supported prior probability distribution over $\xi_j \in (-\infty, \infty)$ with location $\log \tau$ and scale $\psi$. Then, under hierarchy (\ref{eq:GL:Hierarchy:2}) with fixed $y_j$
\begin{equation}
	\label{eq:Kappa:Ineq:1}
	\mathbb{P}(\kappa_j \leq \varepsilon \vbar y_j, \tau, \psi) \leq \left( \frac{\int_{\xi(\varepsilon)}^{\infty} p(\xi \vbar \log \tau, \psi) d \xi}{1 - \int_{\xi(\varepsilon)}^{\infty} p(\xi \vbar \log \tau, \psi) d \xi} \right) e^{y_j^2/2}
\end{equation}
and with fixed $\tau$, $\psi$
\begin{equation}
\label{eq:Kappa:Ineq:2}
	\mathbb{P}(\kappa_j \geq \varepsilon \vbar y_j, \tau, \psi) \leq c \, e^{-y_j^2 \varepsilon (1 - 1/\delta)/2}
\end{equation}
where $c > 0$, $\delta > 1$ and $\xi(\kappa) = (1/2) \log ((1-\kappa)/\kappa)$ is the transformation between $\kappa$ and $\xi$.
\end{proposition}

Inequality (\ref{eq:Kappa:Ineq:2}) states that for large values of $y_j$, most of the posterior probability concentrates away from $\kappa=1$ at an exponential rate, with the obviously corollary that $\mathbb{P}(\kappa \leq \varepsilon \vbar y_j, \tau, \psi) \to 1$ as $|y_j| \to \infty$. This is essentially the same result as Theorem 3.2 in~\cite{DattaGhosh13} which was developed specifically for the horseshoe prior, and demonstrates that performing ``no shrinkage'' is actively controlled by the magnitude of the data for most sensible prior distributions. In contrast, shrinkage to zero is much less determined by $y_j$ being small; rather it is the behaviour of the prior density that determines how aggressively the effects are shrunk towards zero. The implications of inequality (\ref{eq:Kappa:Ineq:2}) are actually quite weak, in the sense that it is valid for any location-scale prior distribution over $\xi_j$ which is supported on $(-\infty,\infty)$; this includes for example, the Bayesian lasso, which has very light tails as $|\beta_j| \to \infty$ and is asymptotically biased for large values of $y_j$. 

Inequality (\ref{eq:Kappa:Ineq:1}) characterises the behaviour of the posterior distribution near $\kappa_j = 0$, i.e., when no shrinkage is applied; an important implication of this inequality is that if $\xi_j$ follows a unimodal location-scale prior distribution with location $\log \tau$, the following limit holds:
\begin{equation}
	\label{eq:Kappa:Ineq:3}
	\mathbb{P}(\kappa_j \leq \varepsilon \vbar y_j, \tau, \psi) \to 0 \mbox{ as } \tau \to 0
\end{equation}
The proof of this follows from the fact that $\int_{a}^{\infty} p(\xi \vbar \log \tau, \psi) d \xi$ is a decreasing function of $\tau$. Result (\ref{eq:Kappa:Ineq:3}) says that as the global shrinkage parameter tends to zero the posterior probability in $\kappa$-space concentrates around $\kappa=1$, i.e., complete shrinkage. That is, for small enough $\tau$, all coefficients can be shrunk arbitrarily close to zero. 

We note that this result is very general, making relatively weak assumptions about the prior distribution $p(\xi_j \vbar \log \tau, \psi)$ and requiring only that we know the CDF of the prior distribution, which is frequently the case. It is therefore of some interest to evaluate inequality (\ref{eq:Kappa:Ineq:1}) for some specific choices of prior distributions. The behaviour of (\ref{eq:Kappa:Ineq:1}) as $\tau \to 0$ in three interesting cases is formalised in the following result.

\begin{proposition}
\label{prop:tau:decreasing}
Under hierarchy (\ref{eq:GL:Hierarchy:2}), with fixed $y_j$:
\begin{enumerate}
	\item if $\lambda_j^2 \sim {\rm Exp}(\tau^2)$, then	$\mathbb{P}(\kappa_j \leq \varepsilon \vbar y_j, \tau, \psi) = O \left( e^{y_j^2/2} \exp \left[ -\frac{(1-\varepsilon)}{\varepsilon \tau^2} \right] \right)$;

	\item if $\lambda_j \sim {\rm LogHS}(\tau, \psi)$, then $\mathbb{P}(\kappa_j \leq \varepsilon \vbar y_j, \tau, \psi) = O \left( e^{y_j^2/2}  \left( \frac{2}{\pi} \right)  \left( \frac{\sqrt{\varepsilon}\tau}{\sqrt{1-\varepsilon}} \right)^{1/\psi} \right)$;

	\item if $\lambda_j \sim \mbox{{\rm Log}-$t$}(\alpha,\tau, \psi)$, then	$\mathbb{P}(\kappa_j \leq \varepsilon \vbar y_j, \tau, \psi) = O \left( e^{y_j^2/2} \left[ \frac{\psi}{ \log \left( \frac{ 1-\varepsilon }{ \varepsilon \tau^2 } \right) } \right]^{\alpha} \right)$;
	
\end{enumerate}
as $\tau \to 0$, where ${\rm LogHS}(\tau, \psi)$ denotes the log-hyperbolic secant distribution (\ref{eq:hyp_sech}) with scale $\tau$ and log-scale $\psi$, and $\mbox{{\rm Log}-$t$}(\alpha, \tau, \psi)$ denotes the log-$t$ distribution (\ref{eq:t:lambda}) with scale $\tau$, log-scale $\psi$ and degree-of-freedom $\alpha$.
\end{proposition}

From proposition (\ref{prop:tau:decreasing}) we see that the upper-bound on the probability of $\kappa$ being close to zero (i.e., no shrinkage) decreases at a different rate depending on the choice of prior for the local shrinkage parameters $\lambda_j$. In the case of the Bayesian lasso, the results suggest that the probability tends to zero much more rapidly than in the case of the log-hyperbolic secant prior (which includes the horseshoe as a special case), which in turn tends to zero much quicker than in the case of the log-$t$ prior with degree-of-freedom $\alpha=1$ (i.e., the log-Cauchy distribution). Use of Proposition \ref{prop:log:concave:bounds} shows that the log-Laplace prior behaves the same as the log-hyperbolic secant prior. While these upper-bounds are not particularly sharp, these results suggest that the log-$t$ prior distribution for $\lambda_j$ yields a posterior distribution over $\kappa_j$ that is much less sensitive to the choice of $\tau$ than the Bayesian lasso or horseshoe. 


The final concentration inequality we will present examines the behaviour of the posterior probability distribution over $\kappa_j$ as the scale parameter $\psi \to 0$. 

\begin{proposition}
\label{prop:Kappa:Ineq:4}
Let $p(\xi_j \vbar \log \tau, \psi)$ be a symmetric, unimodal, fully-supported prior probability distribution over $\xi_j = \log((1-\kappa)/\kappa)/2 \in (-\infty, \infty)$ with location $\log \tau$ and scale $\psi$. Then, under hierarchy (\ref{eq:GL:Hierarchy:2}) with fixed $\tau$ and $y_j$
\begin{equation}
	\label{eq:Kappa:Ineq:4}
	\mathbb{P} \left( \kappa_j \notin ((1+\tau^2 e^{\varepsilon})^{-1}, \, (1+\tau^2 e^{-\varepsilon})^{-1}) \vbar y_j,  \tau, \psi \right) \to 0 
\end{equation}
as $\psi \to 0$, for all $\varepsilon > 0$.

\end{proposition}

Inequality (\ref{eq:Kappa:Ineq:4}) demonstrates that for all symmetric, unimodal prior distributions over $\xi_j$, the posterior probability of the shrinkage coefficient $\kappa_j$ can be concentrated as greatly around $1/(1+\tau^2)$ as desired by making the scale parameter $\psi$ sufficiently close to zero. A corollary of this result is that for small values of the scale parameter $\psi$ the resulting Bayesian procedure essentially mimics ridge regression. To see this, we note that ridge regression can be framed in the local-global hierarchy by taking
%
\[
	\xi_j \sim \delta_{\log \tau}(\xi_j) d\xi_j,
\]
where $\delta_v(x)$ denotes the Dirac point-mass at $x=v$. In this case, $\E{}{\kappa_j \vbar y_j, \tau} = 1/(1+\tau^2)$ which does not depend on the particular value $y_j$, so that the same degree of shrinkage is applied uniformly to all coefficients. 

\subsection{Estimation of $\psi$}
\label{sec:Estimating:psi}

The discussion and results in Sections (\ref{sec:LogScale}) and (\ref{sec:TheoreticalResults}) suggest that introducing a scale parameter $\psi$ into a shrinkage density $f(\xi_j)$ over $\xi_j$ provides a simple, unified method for controlling the concentration and tail behaviour of the resulting marginal distribution over the coefficient $\beta_j$. By making the scale parameter suitably small we can concentrate probability mass near $\kappa_j=1/2$ and obtain near-ridge regression like behaviour. Conversely, making the scale parameter sufficiently large will spread the probability mass over the $\xi_j$ space more thinly, and allow for greater variation in the shrinkage coefficients. To build a shrinkage prior that has the ability to adapt to the degree of sparsity in the underlying coefficient vector we could put an appropriate prior over $\psi$ and incorporate into the global-local prior hierarchy to allow for its estimation along with the other hyperparameters. An advantage of this approach, in comparison to earlier attempts to adaptively control tails by varying shape parameters, is that $\psi$ has a the same interpretation as a scale parameter irrespective of the shrinkage prior $f(\xi_j)$ we begin with. A possible prior for $\psi$ might be
\[
	\psi \sim C^{+}(0,1)
\]
though there exists a large variety of priors for scale parameters in the literature that we can draw upon (see for example, \cite{Gelman06}). 

In the case that $f(\xi_j)$ is log-concave there exists a potential problem when varying the scale $\psi$. Theorem \ref{eq:LogLaplaceTheorem1} suggests that if $f(\xi_j)$ is log-concave, and the scale parameter is too small, the resulting marginal distribution $p(\beta_j)$ can lose Properties I and II discussed in Section \ref{sec:intro:properties}. In fact, the following result shows that for all log-concave densities this is always a possibility.

\begin{proposition}
\label{prop:log:concave:scale}
	If $f(\xi_j)$ is a log-concave density over $\mathbb{R}^+$ with a maxima at $\xi_j = 0$, and $f(\xi_j \vbar \psi)$ is the scale-density
	\[
		f(\xi_j \vbar \psi) = \left( \frac{1}{\psi} \right) f\left( \frac{\xi_j}{\psi} \right)
	\]
	then there always exists a $\psi>0$ such that
	\[
		\sup_{\xi_j < \mu} |g(\xi_j \vbar \psi)| > 1 \; \; \mbox{and} \; \;	\sup_{\xi_j > \mu} |g(\xi_j \vbar \psi)| > 1
	\]
	where $g = - d \log f(\xi_j \vbar \psi)/d \xi_j$.
\end{proposition}

This result shows that if the density $f(\xi_j)$ we use to create our scale-density is log-concave, then there will always exist a choice of scale parameter $\psi>0$ such that the gradients of the log-density exceed one on either side of the mode. We can use this in conjunction with Proposition \ref{prop:log:concave:bounds} to show that we can always find a $\psi^\prime$ such that $f(\xi_j \vbar \psi^\prime)$ is upper-bounded by a log-Laplace density $p_{\rm LL}(\xi_j \vbar \psi_1^\prime, \psi_2^\prime)$ with $\psi_1^\prime < 1$ and $\psi_2^\prime < 1$. Then, application of Theorem \ref{eq:LogLaplaceTheorem1} shows the corresponding marginal density $p(\beta)$ will not have a pole at $\beta_j = 0$, and application of Theorem \ref{eq:LogLaplaceTheorem2} shows that the tails of $p(\beta_j)$ will be sub-Cauchy. The implication of this result is that if we allow estimation of the scale parameter $\psi$ from data, we can only obtain ridge-like behaviour at the cost of losing the ability to estimate large signals without bias and to aggressively shrink away small signals. We can also use the properties of the log-Laplace prior to establish a sufficient condition for a prior distribution over $\xi_j$ to lead to a marginal distribution over $\beta$ that possesses Properties I and II; in particular, from Proposition \ref{prop:log:linear:tails} and Theorems \ref{eq:LogLaplaceTheorem1} and \ref{eq:LogLaplaceTheorem2} we know that it if we can find an asymmetric Laplace distribution that lower-bounds $f(\xi_j)$ and has scale parameters $\psi_1 \geq 1$, $\psi_2 \geq 1$, then the resulting marginal distribution over $\beta_j$ will possess Properties I and II.

If we move outside of the class of log-concave scale prior distributions for $\xi_j$ we see that the results can be strikingly different. For example, if we consider prior distributions for $\lambda_j$ built from the log-$t$ distribution then Theorem \ref{eq:Log:student:t:theorem} shows that these priors always result in marginal distributions $p_t(\beta_j \vbar \psi, \alpha)$ that possess properties I and II, irrespective of the value of $\psi$. Thus, the log-$t$ priors may place as much mass around $\xi_j = 0$, or conversely $\kappa_j = 1/2$, as desired, to restrict the variation in shrinkage coefficients while still providing the possibility of either heavily shrinking coefficients close to zero, or leaving very large coefficients virtually unshrunk. The log-$t$ distribution would therefore appear to offer a family of shrinkage priors that smoothly transition from ultra-sparse to ultra-dense prior beliefs, while being safe in the sense that they always provide an ``out'' for very large coefficients, or coefficients that are exactly zero. This suggests that this class of prior distributions is potentially a strong candidate around which to try and build robust, adaptive shrinkage estimators. To the best of our knowledge, these are the only priors constructed so far that have this particular property.

\subsubsection{Estimation of $\tau$}

A final point of interest regarding the simulations was the observation that if the global shrinkage hyperparameter $\tau^2$ was not constrained to the interval $(1/n,1)$, the adaptive log-$t$ prior almost always estimated $\psi$ to be very small, and had subsequently performed almost identically to ridge regression. This behaviour appears to be related to the well-known degeneracy of the likelihood of the $t$-distribution~\citep{FernandezSteel99a} when there are too few degrees-of-freedom in the ``data points'' (i.e., in our setting, the log-shrinkage parameters $\xi_j = \log \lambda_j$); when $\tau$ is large, the shrinkage coefficients tend to concentrate near $\kappa_i = 0$, as shrinkage towards zero is a passive effect induced by the prior distribution, while a lack of shrinkage is actively controlled by the size of the data $y_j$ (as demonstrated by Proposition~\ref{prop:kappa:ineq}). This means that for large $\tau$, most of the shrinkage coefficients $\kappa_j$ will cluster around zero, and therefore the log-shrinkage parameters $\xi_j = \log \lambda_j$ will be very similar to each other; the degeneracy of the $t$-distribution likelihood will then cause a collapse of $\psi$ towards zero. Restricting $\tau$ to be no larger than one appears to avoid this degeneracy, and provides reasonable inferences regarding $\psi$. Further characterisation of this result is the current focus of ongoing research.

%
%

%


\section{Posterior Computation}
\label{sec:posterior:comp}

Ideally, in addition to possessing favourable theoretical properties, a prior distribution should also result in a posterior distribution from which efficient simulation is possible. An interesting aspect of working in terms of $\xi_j = \log \lambda_j$, rather than directly in terms of $\lambda_j$, is that the conditional distributions $p(\xi_j \vbar \cdots)$ are frequently both unimodal and log-concave, which allows for the use of simple and efficient rejection samplers. To see this, we note that for the global-local shrinkage hierarchy (\ref{eq:Global})--(\ref{eq:Local}) the conditional distribution of $\xi_j$ can be written in the form
\begin{equation}
	\label{eq:Conditional_tj_general}
	p(\xi_j \vbar \cdots) \propto \exp \left( {-m_j e^{-2 \xi_j} - \xi_j} \right) p(\xi_j)
\end{equation}
where $m_j = \beta_j^2/\tau^2/\sigma^2$. The first term in (\ref{eq:Conditional_tj_general}) is log-concave in $\xi_j$, and provided that the prior $p(\xi_j)$ is log-concave, the conditional distribution will also be log-concave. As discussed in Section \ref{sec:Comparison:Standard:Priors} this condition is satisfied by many of the standard shrinkage priors (i.e., horseshoe, horseshoe+, Bayesian lasso). However, even in the case of non-log-concave distributions, such as the Student-$t$, many symmetric location-scale distributions can be expressed as a scale-mixture of normal distributions, i.e.,
\[
	\xi_j \vbar \omega_j^2 \sim N(0, \psi^2 \omega_j^2), \; \; \; \omega_j^2 \sim p(\omega^2) d\omega^2
\]
where $p(\omega_j^2)$ is a suitable mixing density. The normal distribution is log-concave, so that any distribution $p(\xi_j)$ that can be expressed as a scale-mixture of normals admits log-concave conditional distributions for $\xi_j$, conditional on the latent variable $\omega_j$. A further advantage of this representation is that the scale parameter $\psi$, which controls the tails of the induced prior distribution over $\lambda_j$, appears as a simple scale parameter. This means that sampling $\psi$, and adaptively controlling the tail weight of the prior distribution over $\lambda_j$ as discussed in Section \ref{sec:Estimating:psi}, becomes straightforward. 

\subsection{Log-Scale Prior Hierarchy}

In this section we present the steps required to implement a Gibbs sampler for both the symmetric log-Laplace prior density (\ref{eq:LogLaplacePDF}) and the log-$t$ prior density (\ref{eq:t:lambda}), both with unknown scale parameters $\psi$. We use the fact that both of these densities can be written as scale-mixtures of normals in the $\xi_j$ space. We present the hierarchy for the multiple means model (\ref{eq:multiple:means}) with known noise variance $\sigma^2$, though adaptation to the general linear regression model with unknown variance is straightforward as the conditional distributions for the hyperparameters remain the same (see for example \cite{MakalicSchmidt16}). Our hierarchy is
\begin{eqnarray}
	\nonumber
	y_j \vbar \beta_j &\sim& N(\beta_j, \sigma^2) \\
	\nonumber
	\beta_j \vbar \lambda_j, \tau &\sim& N(0, \lambda_j^2 \tau^2) \\
	\nonumber
	\log \lambda_j \vbar \omega_j, \psi &\sim& N(0, \omega_j^2 \psi^2) \\
	\label{eq:Hierarchy}
	\omega_j^2 &\sim& p(\omega_j^2) d \omega_j^2 \\
	\nonumber
	\psi &\sim& C^{+}(0,1) \\
	\nonumber
	\tau &\sim& C^+_{(0, 1)}(0,1)
\end{eqnarray}
where $j = 1,\ldots,n$, and $C^+_{(0,1)}(0,1)$ denotes a standard half-Cauchy distribution truncated to the interval $(0,1)$, as recommended by \cite{VanDerPasEtAl17a}.
The choice of density for the latent variables $\omega_j$ determines the particular log-scale prior to be used. For the log-Laplace prior we use
\[
	\omega_j^2 \sim {\rm Exp}(1)
\]
and for the log-$t$ prior we use
\[
	\omega_j^2 \sim {\rm IG} \left( \frac{\alpha}{2}, \; \frac{\alpha}{2} \right)
\]
where $\alpha$ is the degrees-of-freedom parameter, which we assume to be known. The $z$-distribution can also be represented as a scale mixture, which would allow us to easily extend our hierarchy to an adaptive horseshoe procedure~\citep{BarndorffNielsenEtAl82}. In practice, we use the following inverse gamma-inverse gamma mixture representation of the half-Cauchy prior distribution~\citep{MakalicSchmidt16a}
\[
	\psi^2 \vbar \phi \sim {\rm IG}(1/2, 1/\phi), \; \; \; \phi \sim {\rm IG}(1/2, 1).
\]
This latent variable representation leads to simpler conditional distributions in a Gibbs sampling implementation than the alternative gamma-gamma representation commonly used~\cite{ArmaganEtAl11b}.

\subsection{Gibbs Sampling Procedure}

Given observations ${\bf y} = (y_1,\ldots,y_n)$, we can sample from the posterior distribution using the following Gibbs sampling procedure. The coefficients $\beta_j$, $j=1,\ldots,n$ can be sampled from
	\[
		\beta_j \vbar y_j, \lambda_j, \tau \sim N\left( (1 - \kappa_j) y_j, \; 1 - \kappa_j \right)
	\]
	where $\kappa_j = 1/(1+\lambda_j^2 \tau^2)$. The $\xi_j = \log \lambda_j$ hyperparameters can be sampled using the rejection sampler presented in Appendix II. This sampler is highly efficient, requiring approximately $1.2$ draws per accepted sample in the worst case setting. The $\omega_j^2$ latent variables are sampled according to the particular log-scale prior distribution we have chosen. For the log-Laplace prior $1/\omega_j^2$ is conditionally distributed as per
		\[
			\omega_j^{-2} \vbar \xi_j, \psi \sim {\rm IGauss} \left( \left( {2 \psi^2}/{\xi_j^2} \right)^{1/2}, \; 1/2 \right)
		\]
		where ${\rm IGauss}(\mu,\lambda)$ denotes an inverse Gaussian distribution with mean $\mu$ and shape $\lambda$, while for the log-$t$ prior distribution we sample $\omega_j^2$ from
		\[
			\omega_j^2 \vbar \xi_j, \psi \sim {\rm IG} \left(\frac{\alpha+1}{2}, \; \frac{\xi_j^2}{2 \psi^2} + \frac{\alpha}{2} \right).
		\]
	To sample the scale parameter $\psi$ we use the following conditional densities:
	\begin{eqnarray*}
		\psi^2 \vbar \xi_1,\ldots,\xi_n, \omega_1,\ldots,\omega_n, \phi &\sim& {\rm IG} \left( \frac{n+1}{2}, \; \frac{1}{\phi} + \frac{1}{2} \sum_{j=1}^n \frac{\xi_j^2}{\omega_j^2} \right), \\
		\phi \vbar \psi &\sim& {\rm IG} \left( 1, 1 + \frac{1}{\psi^2} \right).
	\end{eqnarray*}
	Finally, the global shrinkage parameter $\tau$ may be sampled from the truncated conditional posterior using a rejection sampler.	

\section{Experimental Results}

\begin{table*}[t]
\scriptsize
\begin{center}
\begin{tabular}{cccccccccccccc} \toprule[1pt]
~~$q_n$~~ & ~~$A$~~ &  & ~~IGG~~ & ~~HS~~ & ~~HS+~~ & ~~Log-$t$~~ & ~~RR~~ & & ~~IGG~~ & ~~HS~~ & ~~HS+~~ & ~~Log-$t$~~ & ~~RR~~ \\
\multicolumn{3}{c}{} & \multicolumn{5}{c}{Squared Error} & & \multicolumn{5}{c}{Classification Accuracy} \\
\cmidrule{1-14} 
\\
0$\cdot$05 & 2 & & 15$\cdot$70 & 16$\cdot$49 & 16$\cdot$88 & 16$\cdot$79 & 17$\cdot$14 &  & 95$\cdot$09 & 94$\cdot$97 & 95$\cdot$01 & 94$\cdot$65 & 94$\cdot$29 \\
& 4 & & 17$\cdot$88 & 18$\cdot$21 & 17$\cdot$46 & 17$\cdot$50 & 44$\cdot$36 &  & 98$\cdot$51 & 98$\cdot$42 & 98$\cdot$46 & 98$\cdot$75 & 86$\cdot$59 \\
& 8 & & 7$\cdot$89 & 10$\cdot$80 & 9$\cdot$58 & 7$\cdot$30 & 76$\cdot$27 &  & 99$\cdot$49 & 99$\cdot$14 & 99$\cdot$22 & 99$\cdot$74 & 6$\cdot$17 \\
& 16 & & 7$\cdot$37 & 10$\cdot$36 & 9$\cdot$35 & 6$\cdot$50 & 92$\cdot$72 &  & 99$\cdot$52 & 99$\cdot$09 & 99$\cdot$18 & 99$\cdot$80 & 5$\cdot$65 \\
\\
0$\cdot$20 & 2 & & 55$\cdot$71 & 48$\cdot$30 & 51$\cdot$84 & 44$\cdot$15 & 44$\cdot$75 &  & 81$\cdot$88 & 82$\cdot$74 & 82$\cdot$52 & 79$\cdot$08 & 74$\cdot$74 \\
& 4 & & 66$\cdot$49 & 52$\cdot$26 & 49$\cdot$35 & 52$\cdot$15 & 76$\cdot$73 &  & 95$\cdot$21 & 95$\cdot$57 & 96$\cdot$80 & 95$\cdot$86 & 20$\cdot$99 \\
& 8 & & 25$\cdot$35 & 41$\cdot$02 & 35$\cdot$17 & 31$\cdot$04 & 93$\cdot$17 &  & 99$\cdot$60 & 94$\cdot$52 & 97$\cdot$22 & 98$\cdot$30 & 20$\cdot$53 \\
& 16 & & 22$\cdot$57 & 38$\cdot$64 & 34$\cdot$29 & 25$\cdot$62 & 98$\cdot$34 &  & 99$\cdot$57 & 94$\cdot$25 & 96$\cdot$65 & 99$\cdot$00 & 20$\cdot$54 \\
\\
0$\cdot$40 & 2 & & 108$\cdot$44 & 75$\cdot$11 & 83$\cdot$45 & 65$\cdot$20 & 61$\cdot$28 &  & 64$\cdot$35 & 71$\cdot$82 & 69$\cdot$71 & 71$\cdot$65 & 41$\cdot$58 \\
& 4 & & 128$\cdot$96 & 81$\cdot$21 & 77$\cdot$84 & 79$\cdot$27 & 85$\cdot$76 &  & 91$\cdot$10 & 94$\cdot$54 & 95$\cdot$71 & 96$\cdot$34 & 40$\cdot$48 \\
& 8 & & 46$\cdot$42 & 59$\cdot$01 & 55$\cdot$63 & 50$\cdot$67 & 95$\cdot$78 &  & 99$\cdot$69 & 94$\cdot$90 & 96$\cdot$29 & 98$\cdot$24 & 40$\cdot$42 \\
& 16 & & 43$\cdot$00 & 55$\cdot$73 & 53$\cdot$69 & 47$\cdot$17 & 99$\cdot$12 &  & 99$\cdot$70 & 94$\cdot$99 & 96$\cdot$29 & 98$\cdot$64 & 40$\cdot$38 \\
\\
0$\cdot$60 & 2 & & 162$\cdot$12 & 96$\cdot$20 & 109$\cdot$57 & 87$\cdot$95 & 70$\cdot$99 &  & 46$\cdot$69 & 66$\cdot$68 & 60$\cdot$46 & 65$\cdot$47 & 60$\cdot$53 \\
& 4 & & 191$\cdot$05 & 108$\cdot$96 & 104$\cdot$79 & 106$\cdot$86 & 89$\cdot$45 &  & 87$\cdot$12 & 95$\cdot$66 & 96$\cdot$28 & 96$\cdot$66 & 60$\cdot$29 \\
& 8 & & 69$\cdot$79 & 77$\cdot$98 & 74$\cdot$26 & 70$\cdot$15 & 98$\cdot$02 &  & 99$\cdot$79 & 96$\cdot$54 & 97$\cdot$56 & 98$\cdot$90 & 60$\cdot$24 \\
& 16 & & 64$\cdot$02 & 72$\cdot$66 & 71$\cdot$06 & 66$\cdot$25 & 100$\cdot$67 &  & 99$\cdot$80 & 96$\cdot$64 & 97$\cdot$49 & 99$\cdot$11 & 60$\cdot$23 \\
\\\bottomrule[1pt]
\end{tabular}
\caption{Squared-error as a percentage of the least-squares estimator, and percentage of coordinates correctly identified as non-zero or zero, for the inverse-gamma gamma (IGG), horseshoe (HS), horseshoe+ (HS+), adaptive log-$t$ with $\alpha=7$ (Log-$t$) and ridge regression (RR) priors, as the fraction of non-zero elements $q_n$ and signal strength $A$ is varied. For all experiments the total dimensionality of the mean vector is $q=500$.\label{tab:Results1}}
\end{center}
\end{table*}


To investigate how well the adaptive shrinkage procedures based on log-scale priors perform in the multiple means problem we undertook a small simulation study. The basic experimental setup was similar to the simulations undertaken in \cite{BhattacharyaEtAl15} and \cite{BaiGhosh17a}. We characterise the experimental test setting by the fraction of non-zero components, $q_n$, of the of the coefficient vector and the signal to noise ratio. In the recent Bayesian regression literature the focus has been on sparse models with high signal-to-noise ratios, as this is the setting in which sparsity promoting priors achieve the greatest gains over the conventional proportional shrinkage techniques. To test the ability of our proposed log-scale priors to adapt to the configuration of the underlying coefficient vector we specifically expand our simulations to include both dense coefficient vectors as well as low signal-to-noise ratios.

For each choice of the fraction of non-zero components $q_n \in \{ 0.05, 0.2, 0.4 \}$ of $\bm{\beta}$ and signal strength $A = \{2,4,8,16\}$ we generated one hundred data samples from the model (\ref{eq:multiple:means}) with $\sigma^2=1$. For each of these samples we estimated $\bm{\beta}$ using the inverse gamma-gamma (IGG) procedure (\cite{BaiGhosh17a},with $a_n = 1/2 + 1/n$ and $b_n=1/n$), the horseshoe (HS), the horseshoe+ (HS+), the log-$t$ prior (Log-$t$) and the Bayesian ridge (RR). For the log-$t$ prior, we chose $\alpha=7$, as this yields a distribution over $\xi_j$ with the same excess kurtosis as the hyperbolic secant distribution to make the method comparable to the regular horseshoe prior; the scale $\psi$ was given a half-Cauchy prior and sampled along with the other hyperparameters. The squared error between $\bm{\beta}$ and the estimated coefficients, relative to the theoretical squared error of $n$ obtained by the least-squares estimate (i.e., $n$), for each method was averaged over the one hundred tests. The percentage of components of $\bm{\beta}$ correctly classified as either zero, or non-zero, by all priors, using the thresholding rule $\E{}{1-\kappa_j \vbar {\bf y}} > 1/2$, was also computed. We used the posterior mean estimate of $\bm{\beta}$ for the IGG, HS, HS+ and RR, and the posterior median for the log-$t$ prior. This choice was made as the posterior means performed substantially better than the posterior median for the non-adaptive priors, making our comparison as favourable as possible for these procedures.

The results are presented in Table \ref{tab:Results1}. The ridge regression estimator performs the strongest when the underlying coefficient vector is less sparse, and the signal strength is low. In these settings it appears difficult to accurately determine which coefficients to shrink to zero. Ridge regression does not attempt this, and incurs significantly less error when estimating the components of $\bm{\beta}$ that are non-zero than the IGG, HS and HS+. When the signal strength is high or the model is sparse, the Bayesian ridge generally performs much worse than the other four methods. The IGG performs the strongest when the signal strength is high ($A \geq 8$). In this setting it becomes easier to identify which components of $\bm{\beta}$ are exactly zero, and which are non-zero, and apply appropriate levels of shrinkage. Overall, the IGG appears sensitive to the signal strength, and to a lesser degree, the level of sparsity.

The performance of the adaptive log-$t$ prior is interesting. It is the outright winner in only two cases, but demonstrates a high level of robustness to the composition of the underlying coefficient vector. For all settings it achieves squared error and classification accuracies that are very close to those achieved by the best performing method, and it is never ranked worse than the third performing prior in any of the tests. This is in contrast to the other methods, particularly the IGG and ridge, which can perform very poorly, relative to the best performing method, depending on the setting. The HS, ridge and log-$t$ never performed worse than the least-squares estimate, while there are configurations for which HS+ and IGG appear to perform poorer than least-squares. 

When $A=2$, the adaptive log-$t$ procedure generally estimated a scale $\psi$ much smaller than one, and essentially reduced to ridge regression; for $A \geq 8$, the method estimated a $\psi$ substantially larger than one and had behaviour much more similar to the other sparsity promoting priors. The setting of $A=4$ was the most difficult, and the estimated scale was in general very close to one. This represents a compromise between ridge style proportional shrinkage and the more individual shrinkage allowed by the sparsity promoting priors.

\section{Future Work}

A number of extensions to the ideas presented in this paper are immediately apparent.
\begin{itemize}
	\item By varying the choice of the mixing density for the scale-mixture of normal representation of a log-variance prior discussed in Section \ref{sec:posterior:comp}, we can potentially generate a whole range of novel shrinkage priors. For example, if we used $\omega_j \sim C^{+}(0,1)$, the resulting prior distribution $p(\xi_j \vbar \psi)$ would have Cauchy-like tails and be similar in nature to the log-$t$ prior distribution, resulting in extremely heavy tailed prior distributions.
	
	\item Further exploration of the properties of prior distributions over $\xi_j$ that have polynomial tails is required, as these form a novel class of prior distribution that has, to the authors' best knowledge, not yet been carefully examined. 
	
	\item A carefully examination of the prior distribution for the scale parameter $\psi$ for the log-variance priors is required. It is an interesting to ask whether a choice of prior distribution for $\psi$ can be found such that the resulting estimation procedure has excellent performance when the underlying model is sparse while also being minimax in terms of squared-error risk.
	
	\item The extension of these adaptive prior distributions to grouped regression models is of particular interest.	The most basic approach is to associate a scale $\phi_k$ with each of the $K$ groups of predictors to allow for different levels of sparsity within each group. Whether these $K$ scale variables should be given individual prior distributions, or should be tied together through a single prior distribution is an open question. 
\end{itemize}

Both the adaptive log-$t$ and adaptive log-Laplace prior distributions are currently being integrated into the Bayesreg Bayesian regression software packages for MATLAB and R~\citep{MakalicSchmidt16}. 



%
%
%
%

\section*{Appendix I: Proofs}


\noindent {\em Proof of Proposition \ref{prop:log:linear:tails}}. To prove the first statement of the proposition, we note that the assumptions given by (\ref{eq:log:linear:tail:ass:1}) imply that there exist constants $\xi_a<0$, $\xi_b>0$, $K_a>0$ and $K_b>0$ such that
\[
	f(\xi) < K_a (a/2) \exp(a \xi)  \; \; \mbox{for all } \xi<\xi_a, \; \; \mbox{and } f(\xi) < K_b (b/2) \exp(-b \xi) \; \; \mbox{for all } \xi>\xi_b,
\]
where $a>0$, $b>0$. We can then choose the constant
\[
	K = {\rm max}\{K_a,K_b\} \, {\rm max} \left\{ \sup_{\xi \in (\xi_a,\xi_b)} \left\{ \frac{f(\xi)}{p_{\rm LL}(\xi \vbar a, b)} \right\}, 1 \right\}.
\]
which is bounded, as $f(\xi)$ is bounded and $p_{\rm LL}(\xi \vbar a,b)$ is strictly positive; this choice of constant ensures that 
\[
	f(\xi) \leq K \, p_{\rm LL}(\xi \vbar a, b) \; \; \mbox{for all } \xi \in \mathbb{R}.
\]
To prove the first statement, note that assumptions (\ref{eq:log:linear:tail:ass:2}) imply that there exists constants $\xi_a<0$, $\xi_b>0$, $K_a>0$ and $K_b>0$ such that
\[
	K_a \exp(a \xi)/(2 a) < f(\xi) \; \; \mbox{for all } \xi<\xi_a, \; \; \mbox{and } K_b \exp(-b \xi)/(2 b) < f(\xi) \; \; \mbox{for all } \xi>\xi_b.
\]
Similar to above, we can then choose the constant
\[
	K = {\rm min}\{K_a,K_b\} \, {\rm min} \left\{  \sup_{\xi \in (\xi_a,\xi_b)} \left\{ \frac{p_{\rm LL}(\xi \vbar a, b)}{f(\xi)} \right\}^{-1}, 1 \right\}.
\]
which is bounded, as $f(\xi)$ is strictly positive and $p_{\rm LL}(\xi \vbar a,b)$ is bounded; this choice of constant ensures that 
\[
	K \, p_{\rm LL}(\xi \vbar a, b) \leq f(\xi) \; \; \mbox{for all } \xi \in \mathbb{R}
\]
which completes the proof. $\hfill$ $\Box$ \\

\noindent {\em Proof of Proposition \ref{prop:t:dominate:LL}}. First, write the log-$t$ density (\ref{eq:t:lambda}) as $p_t(\lambda_j \vbar \alpha, \psi) = c \, L_t(\lambda_j) / \lambda_j$, where 
\begin{equation}
	\label{eq:L:eq}
	L_t(\lambda_j) = \left( \frac{\log(\lambda_j)^2}{\alpha \psi^2} + 1 \right)^{-(\alpha+1)/2}.
\end{equation}
We first prove that 
\[
	\lim_{\lambda_j \to 0^{+}} \left\{ \frac{p_t(\lambda_j \vbar \alpha, \psi)}{p_{\rm LL}(\lambda_j \vbar \psi_1,\psi_2)} \right\} = \infty.
\]
From (\ref{eq:LogLaplacePDF}) we note that to find the above limit we require only the piece of $p_{\rm LL}(\lambda_j \vbar \psi_1,\psi_2)$ for $\lambda_j<1$. Using this in conjunction with (\ref{eq:t:lambda}) we can write the term inside the limit in the above as
\[	
	\frac{p_t(\lambda_j \vbar \alpha, \psi)}{p_{\rm LL}(\lambda_j \vbar \psi_1,\psi_2)} = \lambda_j^{-1/\psi_1} L_t(\lambda_j)  \; \; \mbox{for } \lambda_j < 1.
\]
Making the substitution $\lambda_j = e^{-z_j}$ we can write the limit as
\[	
	\lim_{z_j \to \infty} \left\{ e^{z_j/\psi_1} \left( 1 + \frac{z^2}{\alpha \psi^2} \right)^{-(\alpha+1)/2} \right\}
\]
which clearly tends to $\infty$ as $z_j \to \infty$. To prove that
%
\[
	\lim_{\lambda_j \to \infty} \left\{ \frac{p_t(\lambda_j \vbar \alpha, \psi)}{p_{\rm LL}(\lambda_j \vbar \psi_1,\psi_2)} \right\} = \infty
\]
we use the piece of $p_{\rm LL}(\lambda_j \vbar \psi_1,\psi_2)$ for $\lambda_j>1$ and write
\[
	\frac{p_t(\lambda_j \vbar \alpha, \psi)}{p_{\rm LL}(\lambda_j \vbar \psi_1,\psi_2)} = \lambda_j^{1/\psi_2} L_t(\lambda_j) \; \; \mbox{for } \lambda_j > 1.
\]
Using the substitution $\lambda_j = e^{z_j}$ we can rewrite the limit as
\[
	\lim_{z_j \to \infty} \left\{ e^{z_j/\psi_2} \left( 1 + \frac{z^2}{\alpha \psi^2} \right)^{-(\alpha+1)/2} \right\}
\]
which clearly tends to $\infty$ as $z_j \to \infty$. $\hfill$ $\Box$ \\


\noindent {\em Proof of Proposition \ref{prop:logt:kappa:concentration}}. First we note that the $\xi_j$ corresponding to shrinkage coefficient $\kappa_j = 1/2 + \epsilon$ is 
\[
	\xi_j(\varepsilon) = \frac{1}{2} \log \left( \frac{1 - 2 \varepsilon}{1 + 2 \varepsilon} \right).
\]
The function $\xi_j(\epsilon)$ is symmetric in the sense that $\xi_j(-\epsilon) = - \xi_j(\epsilon)$; therefore, the we have the following equality:
\begin{equation}
	\label{eq:Integral:kappa:t}
	\int_{1/2-\epsilon}^{1/2+\epsilon} p_t(\kappa_j \vbar \alpha, \psi) d \kappa_j = \int_{-\xi_j(\epsilon)}^{\xi_j(\epsilon)} p_t(\xi_j \vbar \alpha,\psi) d \xi_j,
\end{equation}
where $p_t(\xi_j \vbar \cdot)$ is given by (\ref{eq:t:xi}). As $p_t(\xi_j \vbar \cdots)$ is symmetric and unimodal around $\xi_j = 0$, and $\psi$ is a scale-parameter, we can always find $\psi$ sufficiently small to ensure that the integral (\ref{eq:Integral:kappa:t}) exceeds any choice of $\delta<1$. $\hfill$ $\Box$ \\


\noindent {\em Proof of Proposition \ref{prop:kappa:ineq}}. We make use of the fact that if $x$ is a random variable on $T$, and $S_1,S_2 \subset T$, then $\mathbb{P}(x \in S_1) \leq \mathbb{P}(x \in S_1)/\mathbb{P}(x \in S_2)$. Let $p(\kappa_j \vbar \tau, \psi)$ denote the prior implied by $p(\xi_j \vbar \log \tau, \psi)$ over $\kappa_j = 1/(1+e^{2 \xi_j})$. We first prove (\ref{eq:Kappa:Ineq:1}); using (\ref{eq:posterior:kappa}), and noting that (i) $\kappa_j^{1/2}$ is strictly increasing in $\kappa_j$; (ii) $e^{-\kappa_j y_j^2/2}$ is strictly decreasing in $\kappa_j$, we have
\begin{eqnarray*}
	\mathbb{P}(\kappa_j \leq \varepsilon \vbar y_j, \tau, \psi) &\leq& \frac{ \int_{0}^{\varepsilon} \kappa_j^{1/2} e^{-\kappa_j y_j^2/2} p(\kappa_j \vbar \tau, \psi) d\kappa_j }{\int_{\varepsilon}^{1} \kappa_j^{1/2} e^{-\kappa_j y_j^2/2} p(\kappa_j \vbar \tau, \psi) d\kappa_j } \\
	&\leq& \frac{ \varepsilon^{1/2} \int_{0}^{\varepsilon}  p(\kappa_j \vbar \tau, \psi) d\kappa_j }{\varepsilon^{1/2} e^{-y_j^2/2} \int_{\varepsilon}^{1}  p(\kappa_j \vbar \tau, \psi) d\kappa_j } \\
	&\leq& \left( \frac{ \int_{0}^{\varepsilon} p(\kappa_j \vbar \tau, \psi) d\kappa_j }{\int_{\varepsilon}^{1} p(\kappa_j \vbar \tau, \psi) d\kappa_j } \right) e^{y_j^2/2}.
\end{eqnarray*}
Recalling that $\xi(\kappa) = \log((1-\kappa)/\kappa)/2$ completes the proof. To prove (\ref{eq:Kappa:Ineq:2}) we let $\delta > 1$ be an arbitrary constant, and write
\begin{eqnarray*}
	\mathbb{P}(\kappa_j \geq \varepsilon \vbar y_j, \tau, \psi) &\leq& \frac{ \int_{\varepsilon}^{1} \kappa_j^{1/2} e^{-\kappa_j y_j^2/2} p(\kappa_j \vbar \tau, \psi) d\kappa_j }{\int_{\varepsilon/(2 \delta)}^{\varepsilon/\delta} \kappa_j^{1/2} e^{-\kappa_j y_j^2/2} p(\kappa_j \vbar \tau, \psi) d\kappa_j } \\
	&\leq& \frac{ e^{-y_j^2 \varepsilon/2} \int_{\varepsilon}^{1} p(\kappa_j \vbar \tau, \psi) d\kappa_j }{\left(\frac{\varepsilon}{2 \delta}\right)^{1/2} e^{-y_j^2 \varepsilon /(2 \delta)} \int_{\varepsilon/(2 \delta)}^{\varepsilon/\delta} p(\kappa_j \vbar \tau, \psi) d\kappa_j } \\
	&\leq& \left( \frac{ \int_{\varepsilon}^{1} p(\kappa_j \vbar \tau, \psi) d\kappa_j }{\int_{\varepsilon/(2 \delta)}^{\varepsilon/\delta} p(\kappa_j \vbar \tau, \psi) d\kappa_j } \right) \left( \frac{\varepsilon}{2 \delta} \right)^{1/2} e^{-y_j^2 \varepsilon (1 - 1/\delta)/2}.
\end{eqnarray*}
Noting that the first two terms in the product are independent of $y_j$ completes the proof. $\hfill$ $\Box$ \\


\noindent {\em Proof of Proposition \ref{prop:tau:decreasing}}. To prove the three statements we use the following result. If $p(\xi_j \vbar \log \tau, \psi)$ is a fully-supported unimodal distribution with location $\log \tau$ and scale $\psi$, then $\int_{a}^\infty p(\xi^\prime \vbar \log \tau, \psi) d \xi^\prime \to 0$ as $\tau \to 0$. Using this fact we see that the numerator of (\ref{eq:Kappa:Ineq:1}) tends to zero while the denominator tends to one, and we can write
\begin{equation}
	\label{eq:O:kappa:ineq}
	\mathbb{P}(\kappa_j \leq \varepsilon \vbar y_j, \tau, \psi) = O \left( e^{y_j^2/2} \int_{\xi(\varepsilon)}^\infty p(\xi^\prime \vbar \log \tau, \psi) d \xi^\prime \right).
\end{equation}
where $\xi(\varepsilon) = \log ((1-\varepsilon)/\varepsilon)/2$. To prove statement (1), we note that if $\lambda_j^2 \sim {\rm Exp}(\tau^2)$, then 
\begin{eqnarray*}
	\int_{\xi(\varepsilon)}^\infty p(\xi \vbar \log \tau, \psi) d \xi &=& 
	\exp \left( -\frac{1-\varepsilon}{\varepsilon \tau^2} \right).
\end{eqnarray*}
and using this result in (\ref{eq:O:kappa:ineq}) completes the proof. To prove statement (2), we note that if $\lambda_j \sim {\rm LogHS}(\tau, \psi)$ then from (\ref{eq:hyp_sech}) and the fact that $\lambda = \sqrt{(1-\varepsilon)/\varepsilon}$, along with the series expansion of $\arctan(x)$ (\cite{OldhamEtAl00}, p. 357) we have
\begin{eqnarray*}
	\int_{\xi(\varepsilon)}^\infty p_{\rm SECH}(\xi \vbar \log \tau, \psi) d \xi &=& \left( \frac{2}{\pi \psi \tau} \right) \int_{\sqrt{\frac{1-\varepsilon}{\varepsilon}}}^{\infty} \frac{(\lambda/\tau)^{1/\psi-1} }{\left( (\lambda/\tau)^{2/\psi}+1 \right)} d \lambda, \\
	&=& \left( \frac{2}{\pi} \right) \arctan \left( \left( \frac{\sqrt{\varepsilon} \tau}{\sqrt{1-\varepsilon}} \right)^{1/\psi} \right), \\
	&=& O \left( \left(\frac{2}{\pi}\right) \left( \frac{\sqrt{\varepsilon} \tau}{\sqrt{1-\varepsilon}} \right)^{1/\psi} \right),
\end{eqnarray*}
which combined with (\ref{eq:O:kappa:ineq}) completes the proof.  Finally, to prove statement (3) we note that if $\lambda_j \sim \mbox{{\rm Log}-$t$}(\alpha,\tau, \psi)$ then $\xi_j = \log \lambda_j$ follows a $t$-distribution with degrees-of-freedom $\alpha$, location $\log \tau$ and scale $\psi$. Then, it is well known that
\[
	\int_{\xi(\varepsilon)}^\infty p_t(\xi \vbar \alpha, \log \tau, \psi) d\xi = \frac{1}{2} + \left( \frac{\Gamma((\alpha+1)/2)}{\sqrt{\pi \alpha} \psi \Gamma(\alpha/2)} \right) {}_2 F_1 \left( \frac{1}{2}, \frac{\alpha+1}{2}, \frac{3}{2}, -\frac{(\xi(\varepsilon) - \log \tau)^2}{\alpha \psi^2} \right) (-\xi(\varepsilon) + \log \tau)
\]
where ${}_2F_1(\cdot)$ denotes the Gauss hypergeometric function. Using Mathematica~\citep{} to expand around $\tau = 0$ yields
\begin{eqnarray*}
	\int_{\xi(\varepsilon)}^\infty p_t(\xi \vbar \alpha, \log \tau, \psi) d\xi &=& \left( \frac{\psi}{\xi(\varepsilon) - \log\tau} \right)^\alpha \left( \frac{\alpha^{\alpha/2-1} \Gamma((\alpha+1)/2)}{\sqrt{\pi} \Gamma(\alpha/2)} \right) + O \left( \frac{1}{(-\log \tau)}^{a+2} \right), \\
	&=& O \left( \left( \frac{\psi}{ \log \left( \frac{ 1-\varepsilon }{ \varepsilon \tau^2 } \right) } \right)^{\alpha} \right),
\end{eqnarray*}
%
which completes the proof. $\hfill$ $\Box$ \\


\noindent {\em Proof of Proposition \ref{prop:Kappa:Ineq:4}}. We first note that if $p(\xi_j \vbar \log \tau,\psi)$ is a fully supported, unimodal density with location $\log \tau$ and scale $\psi$ then for all $\varepsilon>0$ and all $\delta \in (0,1)$ there exists a $\psi > 0$ such that
\begin{equation}
	\label{eq:LogScaleConcentration}
	\mathbb{P}(\xi_j \in (\log \tau - \varepsilon, \, \log \tau + \varepsilon)) \geq \delta.
\end{equation}
This result is obvious from the properties of location-scale distributions. We now prove that the posterior probability of $\kappa_j \leq (1+\tau^2 e^{2\varepsilon})^{-1}$ tends to zero as $\psi \to 1$; using inequality (\ref{eq:Kappa:Ineq:1}) from Proposition \ref{prop:tau:decreasing}, and letting $c>0$ denote a constant that does not depend on $\psi$, we can write
\begin{eqnarray*}
	\mathbb{P}(\kappa_j \leq (1+\tau^2 e^{2\varepsilon})^{-1} \vbar y_j, \tau, \psi) &\leq& c \left( \frac{\int_{\log \tau +  \varepsilon}^{\infty} p(\xi \vbar \log \tau, \psi) d\xi}{\int_{-\infty}^{\log \tau + \varepsilon} p(\xi \vbar \log \tau, \psi) d\xi} \right), \\
	&\leq& c \left( \frac{1 - \int_{\log \tau - \varepsilon}^{\log \tau + \varepsilon} p(\xi \vbar \log \tau, \psi) d\xi}{\int_{\log \tau - \varepsilon}^{\log \tau + \varepsilon} p(\xi \vbar \log \tau, \psi) d\xi} \right).
\end{eqnarray*}
Using inequality (\ref{eq:LogScaleConcentration}) we see that the numerator of the last equation in the preceding display tends to zero and the denominator tends to one as $\psi \to 0$. We next prove that the posterior probability of $\kappa_j \geq (1+\tau^2 e^{-2\varepsilon})^{-1}$ also tends to zero as $\psi \to 0$. From (\ref{eq:posterior:kappa}) we have
\begin{eqnarray*}
	\mathbb{P}(\kappa_j \geq (1+\tau^2 e^{-2\varepsilon})^{-1} \vbar y_j, \tau, \psi) &\leq& \left( \frac{e^{-(1+\tau^2e^{-2\varepsilon}) y_j^2/2} }{(1 + \tau^2 e^{2\varepsilon})^{-1/2} e^{-y_j^2/2}}\right) \left( \frac{ \int_{1/(1+\tau^2 e^{-2\varepsilon})}^1 p(\kappa \vbar \tau, \psi) d \kappa }{ \int_{1/(1+\tau^2 e^{2\varepsilon})}^{1} p(\kappa \vbar \tau, \psi) d\kappa} \right), \\
	&\leq& c \left( \frac{ \int_{-\infty}^{\log \tau - \varepsilon} p(\xi \vbar \log \tau, \psi) d \xi }{ \int_{-\infty}^{\log \tau + \varepsilon} p(\xi \vbar \log \tau, \psi) d\xi} \right), \\
	&\leq& c \left( \frac{ 1 - \int_{\log \tau - \varepsilon}^{\log \tau + \varepsilon} p(\xi \vbar \log \tau, \psi) d \xi }{ \int_{\log \tau - \varepsilon}^{\log \tau + \varepsilon} p(\xi \vbar \log \tau, \psi) d\xi} \right),
\end{eqnarray*}
and by using inequality (\ref{eq:LogScaleConcentration}) we see that the numerator of the last equation in the preceding display  tends to zero and the denominator tends to one as $\psi \to 0$, which completes the proof. $\hfill$ $\Box$ \\


\noindent {\em Proof of Proposition \ref{prop:log:concave:scale}}. The derivative of $-\log f(\xi \vbar \psi)$ with respect to $\xi$ is
\begin{eqnarray*}
	-\frac{d}{d \xi} \log f(\xi \vbar \psi) &=& - \frac{d f(\xi/\psi)/d (\xi/\psi)}{\psi f(\xi/\psi)} \\
	&=& - \frac{f^\prime(\phi)}{\psi f(\phi)}
\end{eqnarray*}
where $f^\prime(u) = d f(u)/du$ and $\phi = \xi/\psi$. It follows from the properties of log-concavity and unimodality at $\xi=0$ that $f(\phi) \leq f(0)$ (\cite{SaumardWellner14a}, p. 23) and that there exists both a value of $\phi<0$ and a value of $\phi>0$ such that $f^\prime(\phi) \neq 0$. Letting $\phi^\prime$ denote a value of $\phi < 0$ such that $f^\prime(\phi^\prime) \neq 0$, we have
\[
	\left| \frac{d}{d \xi} \log f(\xi=\phi^\prime \psi \vbar \psi) \right| \geq |f^\prime(\phi^\prime) f(0)^{-1} \psi^{-1}|
\]
which can be made greater than one by taking $\psi < f(0)^{-1} f^\prime(\phi^\prime)$, with a similar argument for $\phi > 0$. $\hfill$ $\Box$ \\


\noindent {\em Proof of Theorem \ref{thm:LLMarginal}}. Using (\ref{eq:LogLaplacePDF}) we can write the marginal density (\ref{eq:MarginalForBeta}) for $\beta_j$ as
\begin{equation}
	\label{eq:A1:A2}
	p_{\rm LL}(\beta_j) = \frac{1}{2 \sqrt{2 \pi}}  \left[ \underbrace{\frac{1}{\psi_1} \int_{0}^{1} \frac{\exp \left( - \beta_j^2/2/\lambda^2_j \right)}{\lambda_j^{2-1/\psi_1}} d \lambda_j}_{A_1(\beta_j \vbar \psi_1)} +  \underbrace{\frac{1}{\psi_2} \int_{1}^{\infty} \frac{\exp \left( - \beta_j^2/2/\lambda^2_j \right)}{\lambda_j^{2+1/\psi_2}} d \lambda_j}_{A_2(\beta_j \vbar \psi_2)} \right]
\end{equation}
We begin by evaluating the term $A_1(\beta_j \vbar \psi_1)$. Making the transformation of variables $\lambda_j = 1/\sqrt{v_j}$ yields
\begin{eqnarray}
	\nonumber
	A_1(\beta_j \vbar \psi_1) &=& \frac{1}{2 \psi_1} \int_{1}^{\infty} \frac{\exp \left( - \beta_j^2 v_j / 2 \right)}{v_j^{(1+\psi_1)/(2\psi_1)}} d v_j, \\
	\label{eq:A1}
	&=& \left( \frac{1}{2 \psi_1} \right) E_{\left(\frac{1+\psi_1}{2 \psi_1}\right)} \left( \frac{\beta_j^2}{2} \right),
\end{eqnarray}
where $E_{n}(x)$ denotes the generalized exponential integral of the form 
\[
	E_{n}(x) = \int_{1}^{\infty} \frac{e^{-x t}}{t^n} d t.
\]
We now evaluate the term $A_2(\beta_j \vbar \psi_2)$. Making the transformation of variables $\lambda_j = (\beta_j^2/2/g_j)^{1/2}$
%
and simplifying yields
\begin{eqnarray}
	\nonumber
	A_2(\beta_j \vbar \psi_2) &=& \frac{1}{2 \psi_2} \int_{0}^{\frac{\beta_j^2}{2}} \left( \frac{2}{\beta_j^2} \right)^{\frac{1+\psi_2}{2 \psi_2}} g_j^{\frac{1+\psi_2}{2 \psi_2} - 1} e^{-g_j} d g_j \\
	\label{eq:A2}
	&=& \left( \frac{1}{2 \psi_2} \right) \left( \frac{2}{\beta_j^2} \right)^{\frac{1+\psi_2}{2 \psi_2}} \gamma \left( {\frac{1+\psi_2}{2 \psi_2}}, \, \frac{\beta_j^2}{2} \right),
\end{eqnarray}
where $\gamma(s,x) = \int_{0}^{x} v^{s-1} e^{-v} dv$ denotes the lower incomplete gamma function. Using these expressions for $A_1$ and $A_2$ in (\ref{eq:A1:A2}) completes the proof. $\hfill$ $\Box$ \\


\noindent {\em Proof of Theorem \ref{eq:LogLaplaceTheorem1}}. To prove the result we first write the marginal density (\ref{eq:MarginalForBeta:LL}) as
\begin{equation}
	\label{eq:p:LL:A1:A2}
	p_{\rm LL} ( \beta_j \vbar \psi_1, \psi_2) = \left( \frac{1}{2 \sqrt{2 \pi}} \right) \left(A_1(\beta_j \vbar \psi_1) + A_2(\beta_j \vbar \psi_2) \right)
\end{equation}
where $A_1(\beta_j \vbar \psi_1)$ and $A_2(\beta_j \vbar \psi_2)$ are given by (\ref{eq:A1}) and (\ref{eq:A2}), respectively. We first show that irrespective of the choice of $\psi_2 > 0$, $A_2(\beta_j \vbar \psi_2) = O(1)$ as $\beta_j \to 0$. Define $x = \beta_j^2/2$, $s_1 = (1+\psi_1)/(2 \psi_1)$ and $s_2 = (1+\psi_2)/(2 \psi_2)$. We can then write
\begin{eqnarray*}
	\lim_{|\beta_j| \to 0} \left\{ A_2(\beta_j \vbar \psi_2) \right\} &=& \lim_{x \to 0} \left\{ A_2(x \vbar \psi_2) \right\}, \\
	&=& \left( \frac{1}{2 \psi_2} \right) \, \lim_{x \to 0} \left\{ \frac{ \gamma(s_2, x) }{x^{s_2}} \right\}, \\
	&=& \left( \frac{1}{2 \psi_2} \right) \, \lim_{x \to 0} \left\{ \frac{ x^{s_2-1} e^{-x} }{s_2 x^{s_2-1}} \right\}, \\
	&=& \frac{1}{1+\psi_2},
\end{eqnarray*}
where the third step follows from L'Hopital's rule and the fact that $d \gamma(s,x)/dx = x^{s-1} e^{-x}$ (\cite{OldhamEtAl00}, p. 467). We therefore have
\[
	p(\beta_j \vbar \psi_1, \psi_2) = \left( \frac{1}{2 \sqrt{2 \pi}} \right) A_1(\beta_j \vbar \psi_1) + O(1)
\]
as $|\beta_j| \to 0$. To prove part 3 of the theorem, we note that if $\psi_1 < 1$, then $s_1 > 1$, and from \cite{GautschiCahill72}, p. 229, we have $E_{s_1}(0) = 2 \psi_1/(1-\psi_1)$, and $A_1(\beta_j \vbar \psi_1) = O(1)$ as $|\beta_j| \to 0$. To prove part 2 of the theorem it suffices to note that if $\psi_1 = 1$, then $s_1 = 1$ and $A_1(\beta_j \vbar \psi_1)$ reduces to
\[
	\frac{1}{2} E \left( \frac{\beta_j^2}{2} \right) = O \left( -\log |\beta_j| \right),
\]
where $E(\cdot)$ denotes the exponential integral~(\cite{GautschiCahill72}, p. 229). Finally, to prove part 1 of the theorem, we first note that if $\psi_1 > 1$, then $s_1<1$, and we can use the relationship between $E_n(\cdot)$ and the upper incomplete gamma function $\Gamma(s,x)$ (\cite{OldhamEtAl00}, p. 381) to write $A_1(\cdot)$ as
\begin{eqnarray*}
	A_1(\beta_j \vbar \psi_1) &=& \left( \frac{1}{2 \psi_1} \right) x^{(s_1-1)} \Gamma(1-s_1, x), \\
														&=& \left( \frac{1}{2 \psi_1} \right) x^{(s_1-1)} \left( \Gamma(1-s_1) - \gamma(1-s_1, x) \right), \\
														&=& \left( \frac{1}{2 \psi_1} \right) x^{(s_1-1)} \left( \Gamma(1-s_1) - \sum_{j=0}^\infty \frac{(-1)^j x^{(1-s_1+j)}}{j! (j+1-s_1)} \right), \\
														&=& \left( \frac{1}{2 \psi_1} \right) \left( x^{(s_1-1)}  \Gamma(1-s_1) - \sum_{j=0}^\infty \frac{(-1)^j x^{j}}{j! (j+1-s_1)} \right), \\
														&=& O \left( |\beta_j|^{-1+1/\psi_1} \right),
\end{eqnarray*}
as $|\beta_j| \to 0$, where the second step follows from the fact that $\Gamma(s) = \Gamma(s,x) + \gamma(s,x)$ and the third step follows from the series expansion of $\gamma(s,x)$ (\cite{OldhamEtAl00}, p. 465). \hfill $\Box$ \\


\noindent {\em Proof of Theorem \ref{eq:LogLaplaceTheorem2}}. We first show that the asymptotic behaviour of the marginal distribution does not depend on $\psi_1$: using the asymptotic expansion of $E_s(x)$ in \cite{OldhamEtAl00}, p. 381, in (\ref{eq:p:LL:A1:A2}) we have
\begin{equation}
	\label{eq:p:LL:thm3:1}
	A_1(\beta_j \vbar \psi_1) = O \left( |\beta_j|^{-2} e^{-\beta_j^2} \right)
\end{equation}
as $|\beta_j| \to \infty$. Define $x = \beta_j^2/2$ and $s_2 = (1+\psi_2)/(2 \psi_2)$. Using $\gamma(s,x) = \Gamma(s) - \Gamma(s,x)$, and the asymptotic expansion for $\Gamma(s,x)$ in (\ref{eq:A2}) yields
\begin{eqnarray}
	\nonumber
	A_2(\beta_j \vbar \psi_2) &=& \left( \frac{1}{2 \sqrt{2 \pi}} \right) x^{-s_2} \gamma(s_2, x), \\
	\nonumber
	&=& \left( \frac{1}{2 \sqrt{2 \pi}} \right) x^{-s_2} (\Gamma(s_2) - \Gamma(s_2,x)),  \\
	\nonumber
	&=& \left( \frac{1}{2 \sqrt{2 \pi}} \right) x^{-s_2} \left( \Gamma(s_2) - x^{s_2-1} e^{-x} \left[ 1 + \frac{s_2-1}{x} + \frac{(s_2-1)(s_2-2)}{x^2} + \cdots \right] \right) \\
	\nonumber
	&=& \left( \frac{1}{2 \sqrt{2 \pi}} \right) x^{-s_2} \Gamma(s_2) - x^{-1} e^{-x} \left[ 1 + \frac{s_2-1}{x} + \frac{(s_2-1)(s_2-2)}{x^2} + \cdots \right] \\
	\label{eq:p:LL:thm3:2}
	&=& O \left( |\beta_j|^{-1 - 1/\psi_2} \right).
\end{eqnarray}
as $|\beta_j| \to \infty$. Using (\ref{eq:p:LL:thm3:1}) and (\ref{eq:p:LL:thm3:2}) in (\ref{eq:p:LL:A1:A2}) completes the proof. \hfill $\Box$ \\

\noindent {\em Proof of Theorem \ref{eq:Log:student:t:theorem}}. To prove the first part of this result, we first show that 
\[
	K \, p_t(\lambda \vbar \alpha, \psi) > p_{\rm LL}(\lambda \vbar \psi_1,\psi_2)
\]
for all $\alpha>0$, $\psi>0$, $\psi_1>0$ and $\psi_2>0$, with $K>0$ a constant not depending on $\lambda$. To prove this, we first note that both $p_t(\lambda \vbar \cdot)$ and $p_{\rm LL}(\lambda \vbar \cdot)$ are bounded, non-negative, continuous functions on the interior of $\mathbb{R}$; these facts, coupled with Proposition \ref{prop:t:dominate:LL} suffice to prove the bounds shown in the previous display. Using the monotone convergence theorem we can then prove that the marginal densities satisfy
\begin{equation}
	\label{eq:p_t:p_LL:bound}
	K \, p_t(\beta \vbar \alpha, \psi) > p_{\rm LL}(\beta \vbar \psi_1, \psi_2)
\end{equation}
which holds for all $\psi>0$, $\alpha>0$, $\psi_1>0$ and $\psi_2>0$. From Theorem \ref{eq:LogLaplaceTheorem1}, we know that $\psi_1>1$ implies that $p_{\rm LL}(\beta \vbar \psi_1,\psi_2) = O \left( |\beta|^{-1+1/\psi_1} \right)$ as $|\beta| \to 0$, and using the bound (\ref{eq:p_t:p_LL:bound}) we therefore have $p_{t}(\beta \vbar \alpha, \psi) = \Omega \left( |\beta|^{-1+1/c} \right)$ for all $c>0$ as was claimed. To prove the second part of the theorem, we can write the log-$t$ density (\ref{eq:t:lambda}) as $p_t(\lambda \vbar \alpha, \psi) = c \, \lambda^{-1} L_t(\lambda)$, where 
\[
	L_t(\lambda) = \left( \frac{\log(\lambda)^2}{\alpha \psi^2} + 1 \right)^{-(\alpha+1)/2}
\]
is a function of slow-variation, in the sense that 
\[
	\lim_{\lambda \to \infty} \left\{ \frac{L_t(t \lambda)}{L_t(\lambda)} \right\} = 1.
\]
Then, by Theorem 6.1 in \cite{BarndorffNielsenEtAl82}, page 157, we have
\begin{equation}
	\label{eq:log:t:Barndorff}
	p_t(\beta \vbar \alpha, \psi) \sim c \, |\beta|^{-1} \left( \frac{\log \left(\beta^2 \right)^2}{\alpha \psi^2} + 1 \right)^{-(\alpha+1)/2}
\end{equation}
where $c$ is a constant independent of $\beta$. The second part of the theorem follows directly from (\ref{eq:log:t:Barndorff}). $\hfill$ $\Box$ \\

\noindent {\em Proof of Theorem \ref{thm:KL:super:eff}}. The three conditions ensure that there exists a $K \in (0,\infty)$ such that
\[
	\pi(\beta) \geq K \, p_{\rm DLL}(\beta \vbar \psi_1, \psi_2)
\]
for all $|\beta| \in (0,\infty)$. For the probability model $y_1,\ldots,y_n \sim N(\beta,1)$ the $\varepsilon$-radius Kullback--Leibler ball centered at $\beta_0 = 0$ is given by $A_{\varepsilon} = \{ \beta : \beta^2/2 \leq \varepsilon \}$; letting $w = 2 \sqrt{\varepsilon}$, the set $A_{\varepsilon}$ is the interval $(-w,w)$. Taking $w \leq 1$, the monotone covergence theorem yields
\begin{eqnarray*}
	\int_{-w}^{w} \pi(\beta) d \beta &\geq& 2 \, \int_{0}^{w} \left( \frac{K}{4 \psi_1} \right) \beta^{-1+1/\psi_1} d \beta \\
	&=& \frac{1}{2} K w^{1/\psi_1}.
\end{eqnarray*}
For a given $w$, using this lower-bound in conjunction with Lemma 1 yields the Ces\`{a}ro risk upper-bound
\begin{eqnarray}
	\nonumber
	R^{\pi}_n(0) &\leq& \frac{w^2}{2} - \frac{1}{n} \log \int_{-w}^{w} \pi(\beta) d\beta \\
	\label{eq:Cesaro:risk:LL}
	&\leq& \frac{w^2}{2} - \frac{1}{n} \log \left( \frac{1}{2} K w^{1/\psi_1} \right).
\end{eqnarray}
Sharpening the bound by optimising for $w$ yields
\[
	w = \left( \frac{1}{n \psi_1} \right)^{1/2},
\]
which is less than one if $n \psi_1 > 1$; plugging this value into (\ref{eq:Cesaro:risk:LL}) completes the proof. $\hfill$ $\Box$ \\


\section*{Appendix II: Implementation Details}

\noindent {\bf Rejection Sampling Algorithm}. For the Bayesian hierarchy (\ref{eq:Hierarchy}) the conditional distribution for $\xi_j$ is proportional to
\begin{equation}
	\label{eq:ConditionalXiDensity}
	p(\xi_j \vbar \cdots) \propto \exp \left( -m_j e^{-2 \xi_j} - \xi_j - \frac{\xi_j^2}{2 v_j} \right) 
\end{equation}
where
\[
	m_j = \frac{\beta_j^2}{2 \tau^2 \sigma^2}, \; \; \; v_j = \omega_j^2 \psi^2.
\]
This density is log-concave and unimodal. We exploit these properties by using a rejection sampler based around a piece-wise proposal distribution made from two exponential distributions on either side of a uniform distribution. The mode of (\ref{eq:ConditionalXiDensity}) is located at
\[
	\xi^\prime = \frac{1}{2} W \left( 4 e^{2 v_j} m_j v_j \right) - v_j,
\]
where $W(\cdot)$ is the Lambert-W (polylogarithm) function. A fast algorithm for $W(\cdot)$ in the range that we require is presented below. Define the functions
\begin{eqnarray*}
	L(\xi) &=& \frac{\xi^2}{2 v_j} + m_j e^{-2 \xi} + \xi, \\
	g(\xi) &=& \frac{\xi}{v_j} - 2 m_j e^{-2 \xi} + 1, \\
	H(\xi) &=& \frac{1}{v_j} + 4 m_j e^{-2 \xi}, \\
\end{eqnarray*}
which are the negative logarithm, the gradient of the negative logarithm and the second derivative of the negative logarithm of (\ref{eq:ConditionalXiDensity}), respectively. We can then define the two points
\[
	x_{0} = \xi^\prime - \frac{0.8}{\sqrt{H(\xi^\prime)}}, \; \; \; \; \; x_{1} = \xi^\prime + \frac{1.1}{\sqrt{H(\xi^\prime)}},
\]
which we use to build our envelope. The break-points for the pieces of our proposal are located at
\[
	\xi_{0} = x_{0} - \frac{L(x_{0}) - L(\xi^\prime)}{g(x_{0})}, \; \; \; \; \; \xi_{1} = x_{1} - \frac{L(x_{1}) - L(\xi^\prime)}{g(x_{1})}.
\]
To the left of $x_{0}$, and to the right of $x_{1}$, our proposal follows two different exponential distributions, and in between $x_{0}$ and $x_{1}$ our proposal is a uniform distribution. As (\ref{eq:ConditionalXiDensity}) is log-concave, our proposal can be made to upper-bound the conditional density $p(\xi_j \vbar \cdots)$ with appropriate scaling~\citep{GilksWild92}. Define the quantities
\[
	K_i = \frac{e^{- L(x_i) + L(\xi^\prime) - g(x_i) (\xi_{i} - x_{i} )}}{|g(x_{i})|}, \; \;\;  i \in \{0,1\},
\]
and $K = K_0 + K_1 + (\xi_1-\xi_0)$. Then, we have the following algorithm:
\begin{enumerate}
	\item First, generate the random variables
	\[
		u_1 \sim U(0,1), \; \; \; u_2 \sim U(0,1), \; \; \; u_3 \sim U(0,1),
	\]
	where $U(a,b)$ denotes a uniform distribution on $(a,b)$.
	
	\item Next, generate from the proposal distribution: if $u_1 < (K_0+K_1)/K$:
	\begin{enumerate}
		\item check if $u_1 < K_0/K$; if so, set $i \leftarrow 0$, else set $i \leftarrow 1$;
		
		\item then set
		\[
			\xi_j \leftarrow -\frac{\log (1-u_2)}{g(x_i)} + \xi_i, \; \; \; \; f \leftarrow L(x_i) + g(x_i) (\xi - x_i);
		\]
	\end{enumerate}
	
	otherwise, we set
	\[
		\xi_j \leftarrow u_2 (\xi_1 - \xi_0) + \xi_0, \; \; \; \; f \leftarrow L(\xi^\prime).
	\]

	\item To determine whether we accept $\xi_j$, check if
	\[
		\log u_3 < f - L(\xi_j);
	\]
	if so, we accept $\xi_j$; otherwise, we reject $\xi_j$ and return to step 1.
\end{enumerate}

For the worse-case values of $m_j>0$ and $v_j>0$ this accept-reject algorithm requires around $1.2$ draws from the proposal per accepted sample and is computational efficient. \\

\noindent {\bf Algorithm to Calculate the Lambert-W function}. We summarise the simple algorithm which is adapted from () with improved initialisation, to efficiently compute the Lambert's W-function. To calculate $W(x)$ for $x>0$:
\begin{enumerate}
	\item If $x < 3$, set $w^{(0)} \leftarrow 1$; else $w \leftarrow \log x - \log \log x$

	\item $t \leftarrow 0$
	
	\item $v \leftarrow w^{(t)} e^{w^{(t)}} - x$

	\item $w^{(t+1)} \leftarrow w^{(t)} - v \left( e^{w^{(t)}+1} - v (w^{(t)}+2) (2 w^{(t)}+2)^{-1} \right)^{-1}$

	\item If $|w^{(t)} - w^{(t+1)}|/|w^{(t+1)}| < \varepsilon$ quit algorithm, else set $t \leftarrow t+1$ and go to step 3.
\end{enumerate} 

\bibliography{bibliography}

\end{document}